\documentclass[12pt]{article}
\usepackage{float}
\usepackage{graphicx,tikz,color,url}
\usepackage{amsmath,amssymb,latexsym}
\usepackage{mathdots}
\usetikzlibrary{arrows}
\makeatletter

\newtheorem{theorem}{Theorem}[section]
\newtheorem{definition}[theorem]{Definition}
\newtheorem{lemma}[theorem]{Lemma}
\newtheorem{proposition}[theorem]{Proposition}

\newtheorem{corollary}[theorem]{Corollary}
\newtheorem{remark}[theorem]{Remark}

\newcommand{\proof}{\noindent{\bf Proof.\ }}
\newcommand{\qed}{\hfill $\square$ \bigskip}
\newcommand{\cp}{\,\square\,}
\newcommand{\gsmb}{\gamma_{\rm SMB}}
\newcommand{\gmb}{\gamma_{\rm MB}}
\newcommand{\smallqed}{{\tiny ($\Box$)}}
\newcommand{\po}{{\cal P}_{\rm o}}
\newcommand{\pco}{{\cal P}_{\rm co}}

\newcommand{\cG}{{\cal G}}
\newcommand{\cH}{{\cal H}}
\newcommand{\cF}{{\cal F}}
\newcommand{\cD}{{\cal D}}
\newcommand{\cS}{{\cal S}}
\newcommand{\cP}{{\cal P}}
\newcommand{\cA}{{\cal A}}
\newcommand{\cB}{{\cal B}}
\newcommand{\cT}{{\cal T}}

\newcommand{\w}{{\rm w}}

\begin{document}
	
	\title{Maker-Breaker domination game on Cartesian products of graphs}
	
	\author {Pakanun Dokyeesun\thanks{Email: \texttt{papakanun@gmail.com
	}} }
	\maketitle
	
	\begin{center}
	Institute of Mathematics, Physics and Mechanics, Ljubljana, Slovenia\\
		\medskip

		%
	\end{center}
	
	\begin{abstract}
		The Maker-Breaker domination game is played on a graph $G$ by two players, called Dominator and Staller.  They alternately select an unplayed vertex in $G$. Dominator wins the game if he forms a dominating set while Staller wins the game if she claims all vertices from a closed neighborhood of a vertex. 
		The game is called \emph{D-game} if Dominator starts the game and it is an \emph{S-game} when Staller starts the game.  If Dominator is the winner in the D-game (or the S-game), then $\gmb(G)$ (or $\gmb'(G)$) is defined by the minimum number of moves of Dominator to win the game under any strategy of Staller.
		Analogously, when Staller is the winner, $\gsmb(G)$ and $\gsmb'(G)$ can be defined in the same way.
		
		We determine the winner of the game on the Cartesian product of paths, stars, and complete bipartite graphs, and how fast the winner wins. We prove that Dominator is the winner on $P_m \square P_n$ in both the D-game and the S-game, and $\gmb(P_m \square P_n)$ and $\gmb'(P_m \square P_n)$ are determined when $m=3$ and $3 \le n \le 5$. Dominator also wins on $G \square H$ in both games if $G$ and $H$ admit nontrivial path covers. Furthermore, we establish the winner in the D-game and the S-game on $K_{m,n} \square K_{m',n'}$ for every positive integers $m, m',n,n'$.
		We prove the exact formulas for $\gmb (G)$, $\gmb'(G)$, $\gsmb(G)$, and $\gsmb'(G)$ where $G$ is a product of stars.
	\end{abstract}
	
	
	\noindent
	{\bf AMS Subj.\ Class.\ (2010)}: 91A24, 05C57, 05C65, 05C69
	
	\section{Introduction}
	For a positive integer $n$, we set $[n] = \{1, \ldots, n\}$. Let $G =(V(G), E(G))$ be a graph. The order of $G$ is denoted by $n(G)$. 
	A graph $H$ is a \emph{subgraph} of $G$ if $V(H) \subseteq V(G)$ and $E(H)\subseteq E(G)$. For any $v \in V(G)$, the \emph{open neighborhood} $N_G(v)$ of $v$ is the set of all vertices adjacent to $v$ and the \emph{closed neighborhood} of $v$ is $N_G[v] = N_G(v) \cup \{v\}$. If $S \subseteq V(G)$, then $N_G(S)=\bigcup_{v \in S} N_G(v)$ and $N_G[S]=\bigcup_{v \in S} N_G[v]$.  A set $D \subseteq V(G)$ is a \emph{dominating set} if $N_G[D]=V(G)$ and the minimum cardinality of dominating sets is the \emph{domination number} $\gamma(G)$ of $G$.	
	A set $M \subseteq E(G)$ is a \emph{matching} if no two edges share a vertex in $M$. A matching $M$ is  a \emph{perfect matching} if $M$ covers every vertex in $G$. 	A \emph{path cover} of $G$ is a set of pairwise vertex-disjoint paths which cover $V(G)$. A path cover of $G$ in which every path has length at least $1$ is called a \emph{nontrivial path cover}.
	
	The \emph{Cartesian product} $G\cp H$ of graphs $G$ and $H$ is defined on the vertex set $V(G)\times V(H)$ such that two vertices $(g,h)$ and $(g',h')$ are adjacent if either $gg'\in E(G)$ and $h=h'$, or $g=g'$ and $hh'\in E(H)$. If $h\in V(H)$, then the subgraph of $G\cp H$ induced by the vertex set $\{(g,h):\ g\in V(G)\}$ is a {\em $G$-layer}, and denoted by $G^h$. Analogously the $H$-layers are defined and denoted by $^g\!H$ for a fixed vertex $g\in V(G)$. 
	
	A \emph{hypergraph} $\cH=(V(\cH), E(\cH))$ consists of the vertex set $V(\cH)$ and the \emph{(hyper)edge} set $E(\cH)$ containing nonempty subsets of $V(\cH)$ of any cardinality. In other words, $E(\cH) \subseteq 2^{V(\cH)}$.  A vertex set $T \subseteq V(\cH)$ is a \emph{transversal} (or vertex cover) in $\cH$ if every (hyper)edge in $\cH$ contains at least one vertex in $T$. Note that a loop less graph is a hypergraph where each (hyper)edge contains exactly two vertices. Many of the fundamental definitions associated with graphs can be extended to hypergraphs (for more details, see~\cite{berge}).
	
	The Maker-Breaker game is a positional game introduced in 1973 by Erd\H{o}s and Selfridge~\cite{erdos-1973} and was widely studied (see~\cite{beck-1981, hefetz-2014}). The game is played on a hypergraph $\cH$ by two players, named Maker and Breaker, where the hyperedges of $\cH$ are the \emph{winning sets}. During the game, the players alternately select unplayed vertices in the hypergraph and Maker wins if he occupies all the vertices of a winning set, otherwise Breaker wins.
	
	In 2020, the \emph{Maker-Breaker domination game} (MBD game) was introduced by Duch\^{e}ne, Gledel, Parreau, and Renault in \cite{duchene-2020}. The game is played on a graph $G$ by two players, called Dominator and Staller. Both players alternately select an unplayed vertex in $G$. Dominator wins the game if he can form a dominating set while Staller wins if she can prevent Dominator to form a dominating set. In other words, Staller wins if she claims a closed neighborhood of a vertex in $G$. The total version of the game was studied in~\cite{forcan-2022, gledel-2020}.
	
	The MBD game can be considered as a variation of the Maker-Breaker game where the winning sets are minimal dominating sets, Dominator is Maker and Staller is Breaker. Conversely, if the closed neighborhoods of the vertices are considered to be the winning sets, Dominator is Breaker and Staller is Maker, see~\cite{bujtas-2021, bujtas-2023}. The names of the players are selected to be consistent with the domination game. For more details and results on the domination game, see~\cite{bresar-2010, book-2021, jacobson-1984, kinnersley-2013}. An MBD game is referred to as \emph{D-game} if Dominator is the first player in the game and the game is called \emph{S-game} when Staller starts the game. The sequence  $d_1, s_1, d_2, s_2, \ldots$ is a sequence of played vertices in a D-game and the sequence  $s'_1, d'_1, s'_2, d'_2, \ldots$ is a sequence of played vertices in an S-game respectively.
	
	In~\cite{gledel-2019},  the \emph{Maker-Breaker domination number} (MBD-number) was introduced in the following way.  
	If Dominator has a winning strategy in the D-game, the MBD-number $\gmb(G)$ represents the minimum number of moves Dominator needs to win the game when both players play optimally. Otherwise, when Dominator cannot win the game, $\gmb(G) = \infty$. Similarly, $\gmb'(G)$ is defined in the same way for the S-game. 
	
	In~\cite{bujtas-2021}, the \emph{Staller-Maker-Breaker domination number} (SMBD-number) was defined analogously to the MBD-number. That is, $\gsmb(G)$ is the minimum number of moves Staller needs to win the D-game if both players play optimally. If Staller does not have a winning strategy in the D-game, we set $\gsmb(G)= \infty$. For the S-game the corresponding invariant is denoted by $\gsmb'(G)$.
	
	Recently, Forcan and Qi \cite{forcan-qi-2022} studied the Maker-Baker domination game on Cartesian products of graphs. It was shown that if Dominator wins on a graph $G$ in the D-game and the S-game, then Dominator also wins on $G \cp H$ in both games for every graph $H$. In particular, Dominator always wins on $P_2 \cp H$ in the D-game and the S-game. This inspired us to consider the winner of the game on $P_3 \cp H$ for any graph $H$. 
	To approach the problem, we use nontrivial path covers to investigate the winner on Cartesian products of graphs.
	
	\paragraph{Structure of the paper.}
	In this paper, we study the Maker-Breaker domination game on Cartesian products of paths and complete bipartite graphs. In the next section, we first provide basic properties of the game and recall results which are used in the rest of the paper.
	Then, in Section $3$, we determine the winner of the D-game and the S-game on $P_m \cp P_n$ and $G \cp H$, where $G$ and $H$ admit nontrivial path covers. Moreover, formulas for  $\gmb(P_3 \cp P_n)$ and  $\gmb'(P_3 \cp P_n)$ are proved for $n \in \{3,4\}$. Not all graphs admit a nontrivial path cover, for example, stars. In Section $4$, the outcome of games on products of complete bipartite graphs is studied and it is determined how fast the winner can win the game on products of stars. 
	
	\section{Preliminaries}
	By the definition of the MBD game exactly one player wins in the D-game (or the S-game). The \emph{outcome} of the MBD game for a graph $G$ is defined in~\cite{duchene-2020} according to the winners of the D-game and the S-game.
	\begin{definition}\label{def: outcome}
		The outcome $o(G)$ of $G$ is one of the following:
		\begin{itemize}
			\item[$(i)$] $\cal D$, if Dominator has a winning strategy in the D-game and the S-game,
			\item[$(ii)$] $\cal S$, if Staller has a winning strategy in the D-game and the S-game,
			\item[$(iii)$] $\cal N$, if the first player has a winning strategy in both games.
		\end{itemize}
	\end{definition}

	These are all possible outcomes of the MBD game because the remaining case is not possible as it was shown in \cite{hefetz-2014}.
	
	The disjoint union of graphs $G$ and $H$ is denoted by $G \cup H$.  The following theorem shows the outcome on $G \cup H$.
	
	\begin{theorem}[\cite{duchene-2020}]
		\label{thm:union}
		Let $G$ and $H$ be graphs. Then
		\begin{itemize}
			\item If $o(G)= \cal S$ or $o(H)= \cal S$, then $o(G \cup H)= \cal S$.
			\item If $o(G)=o(H)= \cal N$, then  $o(G \cup H)= \cal S$.
			\item If $o(G)=o(H)= \cal D$, then  $o(G \cup H)= \cal D$.
			\item Otherwise, $o(G \cup H)= \cal N$.
		\end{itemize}
	\end{theorem}

	\begin{theorem}[\cite{gledel-2019}]
		\label{thm:gmb-union}
		If $G$ and $H$ are graphs, then the following hold.
		\begin{itemize}
			\item $\gmb(G)+\gmb(H) \le \gmb(G \cup H) \le \min\{\gmb'(G)+\gmb(H), \gmb(G)+\gmb'(H)\}$.
			\item $\max\{\gmb'(G)+\gmb(H), \gmb(G)+\gmb'(H)\} \le \gmb'(G \cup H) \le \gmb'(G)+\gmb'(H)$.
		\end{itemize}
	\end{theorem}
	
	Given an edge $e$ of a graph $G$, the graph $G-e$ is the subgraph of $G$ where $e$ has been removed. Note that $N_{G-e}[x] \subseteq N_G[x]$ for every $x \in V(G)$. By Proposition 2.2 in \cite{bujtas-2021}, it implies that deleting $e$ is not a disadvantage for Staller. Thus, if Staller wins the D-game and the S-game on $G$, then she also wins the games in $G-e$. On the other hand, if Dominator wins the D-game and the S-game on $G-e$, then he also wins the games in $G$. This result implies the following lemma which will be used often later.
	
	\begin{lemma}\label{lem:o_edgedelete}
		Let $G$ be a graph and $e \in E(G)$.
		\begin{itemize}
			\item [(i)]  If $o(G-e) = \cal D$ then $o(G)= \cal D$. Moreover, $\gmb(G-e) \ge \gmb(G)$ and $\gmb'(G-e) \ge \gmb'(G)$.
			\item  [(ii)]  If $o(G) = \cal S$ then $o(G-e)= \cal S$. Moreover, $\gsmb(G-e) \le \gsmb(G)$ and $\gsmb'(G-e) \le \gsmb'(G)$.
		\end{itemize}
	\end{lemma}

	The pairing strategy of Maker~\cite{hefetz-2014} ensures that he can win in the Maker-Breaker game if an appropriate set of vertex pairs can be defined. Its immediate consequence for the MBD game shows that Dominator has a winning strategy in both the D-game and S-game if the graph admits a perfect matching. Here we state the lemma in a more general form, name as the generalized pairing strategy.

	\begin{lemma}[\cite{bujtas-2021}]
		\label{lem:pairing}
		Consider an MBD game on $G$ and let $X$ and $Y$ be the sets of vertices played by Dominator and Staller, respectively, until a moment during the game. If there exists a matching $M$ in $G-(X\cup Y)$ such that $V(G) \setminus V(M) \subseteq N_G[X]$,
		then Dominator has a strategy to win the continuation of the game, no matter who plays the next vertex.
	\end{lemma}
	\begin{remark}\label{re:pairing}
		To win the continuation of the game, Dominator applies the following strategy. If Staller claims an unplayed vertex $v$ such that $uv \in E(M)$, Dominator responds by playing vertex $u$ it it is unplayed. Otherwise Dominator plays an arbitrary vertex. By this strategy, the number of moves of Dominator is at most $|X|+|M|$.
	\end{remark}

	\begin{lemma} [No-Skip Lemma \cite{gledel-2019, bujtas-2021}]
		\label{lem:no-skip}
		Let $G$ be a graph.
		\begin{itemize}
			\item[(i)] In an optimal strategy of Dominator to achieve $\gmb(G)$ or $\gmb'(G)$ it is never an advantage for him to skip a move. Moreover, if Staller skips a move it can never disadvantage Dominator.
			\item[(ii)] In an optimal strategy of Staller to achieve $\gsmb(G)$ or $\gsmb'(G)$ it is never an advantage for her to skip a move. Moreover, if Dominator skips a move it can never disadvantage Staller.
		\end{itemize}
	\end{lemma}
	As a further consequence of No-skip Lemma, S-game is the D-game when Dominator skips the first move. Similarly, D-game is the S-game when Staller skips the first move. These facts imply the following consequence of No-Skip Lemma.
	
	\begin{corollary} [\cite{gledel-2019, bujtas-2021}]\label{cor:skip} If $G$ is a graph, then
	$\gmb(G) \le \gmb'(G)$ and $\gsmb'(G) \le \gsmb(G).$

	\end{corollary}
	
	In \cite{gledel-2019}, sharp bounds for $\gmb(G)$ and $\gmb'(G)$ were determined. As each vertex is played at most once during the game,  $\gmb(G) < \infty$ implies
	
	\begin{align}\label{bound1}
		1 \le \gmb(G) \le \left\lceil \frac{n(G)}{2} \right\rceil.
	\end{align}
	Similarly, if $\gmb'(G) < \infty$, then
	\begin{align}\label{bound2}
		1 \le \gmb'(G) \le \left\lfloor \frac{n(G)}{2} \right\rfloor.
	\end{align}

	Forcan and Qi investigated the MBD-number on the Cartesian product of two graphs when Dominator is the winner on at least one of these two graphs in the D-game and the S-game. We can rewrite this result as follows.
	\begin{theorem}[\cite{forcan-qi-2022}]
		\label{thm:cart}
		Let $G$ and $H$ be two graphs. If $o(G) = \cal D$ or $o(H) = \cal D$, then $o(G\cp H) = \cal D$.
	\end{theorem}
	This result implies that $o(P_{2m} \cp P_n)=D$ for every positive integer $m, n$. The main result of~\cite{forcan-qi-2022} asserts the following;
	\begin{theorem}[\cite{forcan-qi-2022}]
		\label{thm:P2}
		$\gmb'(P_2 \cp P_n) = n$ for $n \ge 1$ and $\gmb(P_2 \cp P_n) = n-2$ for $n \ge 13$.
	\end{theorem}
	To further investigate the outcome of the game on Cartesian products of paths and complete bipartite graphs, we will use nontrivial path covers.
	
	\begin{theorem}[\cite{Lavasz-1970}] \label{thm: factor}
		A graph $G$ has a nontrivial path cover if and only if $i(G-S) \le 2 |S|$ for every $S \subseteq V(G)$ where $i(G-S)$ is the number of isolated vertices in $G-S$.
	\end{theorem}

	By Theorem~\ref{thm: factor}, we can conclude that $K_{r,s}$ does not have a nontrivial path cover if $1 \le 2r <s$.
	
	\section{Grid graphs}\label{sec:grid-MBD}
	In this section, we set $Z = P_m\cp P_n$,  where $m, n$ are positive integers, $V(Z) = \{(i,j):\ i\in [m], j\in [n]\}$, and $E(Z) = \{(i,j)(i,j+1):\ i\in [m], j\in [n-1]\} \cup \{(i,j)(i+1,j):\ i\in [m-1], j\in [n]\}$. We consider the outcome of the MBD game on $Z$ and the MBD-number of some small grids. By completing previous partial results, we  prove the following characterization for the outcome of the game on grids.
	
	\begin{theorem} \label{thm:pmpn}
		If $n \ge m \ge 2$, then $o(P_m \cp P_n)= \cal D$.
	\end{theorem}
	
	\proof Assume that $n \ge m \ge 2$ and $Z = P_m \cp P_n$. We will show that Dominator has a winning strategy in the D-game and the S-game on $Z$. By Corollary \ref{cor:skip}\,, it suffices to show that Dominator wins the S-game on $Z$.
	
	\paragraph{Case 1.} $m=2$. Then Dominator wins the S-game by Theorem~\ref{thm:cart}\,.
	
	\paragraph{Case 2.} $m=n=3$. Consider the first move $s'_1$ of Staller. By symmetry, it is enough to consider the following cases.
	\begin{itemize}
		\item  If $s'_1 \in \{(1,1)$, $(2,2)\}$, then Dominator replies by playing $d'_1 =(1,2)$. We can find a matching $M$ in $Z-\{s'_1, (1,2)\}$ such that $V(Z)\setminus V(M) \subseteq N_Z[(1,2)]$. By Lemma~\ref{lem:pairing}\,, Dominator has a winning strategy in the game.
		\item  If $s'_1=(1,2)$, then Dominator replies by playing $d'_1=(2,1)$.
		In the second turn, if $s'_2 \in \{(1,1), (2,2), (2,3), (3,1)\}$, then Dominator plays $d'_2=(1,3)$. One can see that only $(3,2)$ and $(3,3)$ remain undominated after the moves $d'_1$ and $d'_2$. Then Dominator will win the game by playing one of these two vertices in his next move. 
		Suppose that $s'_2 \in \{(1,3), (3,2), (3,3)\}$. Then Dominator replies by playing $d'_2=(2,3)$. Thus only vertices $(1,2)$ and $(3,2)$ remain undominated after the moves $d'_1$ and $d'_2$. If $s'_3 \ne (2,2)$, then Dominator replies $d'_3 = (2,2)$ and he wins the game.
		Assume that $s'_3=(2,2)$. Dominator plays $d'_3 \in \{ (1,1), (1,3)\}$ and then he will play an unplayed vertex in the layer $^3 P_3$ in his next move.
		In all cases, Dominator wins the game within his next two moves.
	\end{itemize}
	Hence, Dominator has a winning strategy in the S-game on $P_3\cp P_3$.
	
	\paragraph{Case 3.} $m \ge 3$ and $n \ge 4$. Observe that each path $P_{\ell}$, where $\ell \ge 3$ admits a path cover obtained from only copies of $P_2$ and $P_3$.
	Let $X$ and $Y$ be path covers of $P_m$ and $P_n$, respectively, such that they consist of copies of $P_2$ and $P_3$.
	Let $P_m'$ and $P_n'$ be the disjoint union\index{disjoint union} of paths from the path covers $X$ and $Y$, respectively. Set $Z'= P_m' \cp P_n'$. Then $Z'$ is a disjoint union of copies of $P_2 \cp P_2$, $P_2 \cp P_3$, $P_3 \cp P_2$, and $P_3 \cp P_3$. By Case 1.\,and Case 2., Dominator wins the S-game on every component of $Z'$. By Theorem~\ref{thm:union}\,, Dominator can win the game on $Z'$. Recall that $Z'$ is a spanning subgraph of $Z = P_m \cp P_n$ that is, $Z'$ can be obtained by repeatedly deleting edges from $Z$. Lemma~\ref{lem:o_edgedelete} then implies that Dominator can win the S-game on $Z$.
	\qed
	\medskip 
	
	According to Theorem~\ref{thm:pmpn}\,, we can conclude the outcome of the MBD game on Cartesian products of graphs which admit nontrivial path covers as follows.
	\begin{theorem}\label{thm:nontrivial-GH}
		If $G$ and $H$  are graphs which admit nontrivial path covers, then $o(G \cp H) = \cal D$.
	\end{theorem}
	
	\proof
	Let $X$ and $Y$ be nontrivial path covers of $G$ and $H$, respectively. Let $G'$ and $H'$ be the disjoint union of paths from the path covers $X$ and $Y$, respectively. Then $G' \cp H'$ is a disjoint union of copies of Cartesian products of nontrivial paths.  By Theorem~\ref{thm:pmpn}\,, Dominator wins the MBD games on every component in $G' \cp H'$. It implies that Dominator wins the games on $G' \cp H'$ by Theorem~\ref{thm:union}\,.
	Observe that $G' \cp H'$ is a spanning subgraph of $G \cp H$. By Lemma~\ref{lem:o_edgedelete}\,, it follows that Dominator wins the games on $G \cp H$.
	\qed
	\medskip 
	
	We now consider exact formulas for MBD-numbers on grids $P_3 \cp P_n$ for every positive integer $n \ge 2$.
	\begin{proposition}\label{prop:p3p3}
		$\gmb (P_3 \cp P_3) = \gmb' (P_3 \cp P_3) =4$.
	\end{proposition}
	\proof Now $Z = P_3 \cp P_3$. By Theorem~\ref{thm:pmpn}\,, we know that Dominator has a winning strategy in the D-game and the S-game.
	By Corollary~\ref{cor:skip} and inequality (\ref{bound2}),  $ \gmb(Z) \le \gmb'(Z) \le \lfloor\frac{9}{2} \rfloor = 4$.
	
	It remains to show that  $\gmb(Z) \ge 4$, that is, Staller has a strategy to ensure that Dominator cannot win the D-game within three moves. By symmetry, it is enough to consider the following three cases.
	\begin{itemize}
		\item If $d_1=(1,1)$, then Staller replies by playing $s_1=(3,2)$. If $d_2 \notin \{(2,1), (3,1)\}$, then Staller responds by choosing $s_2=(3,1)$. If $d_2 \in \{(2,1), (3,1)\}$, then Staller selects $s_2=(2,3)$. In any case, Dominator needs to play at least two more vertices to dominate $Z$.
		
		\item If $d_1=(1,2)$, then Staller replies by playing $s_1=(3,1)$. If $d_2 \in \{(1,3), (2,3), (3,3)\}$, then Staller responds by choosing  $s_2=(2,1)$. If $d_2 \in \{(1,1), (2,1)\}$, then Staller replies by playing  $s_2=(3,3)$. If $d_2 =(3,2)$, then Staller selects $s_2=(2,2)$. Otherwise, Staller plays $s_2=(3,2)$. By this strategy, Dominator needs to plays two more vertices to dominate $Z$.
		
		\item If $d_1=(2,2)$, then Staller replies by playing $s_1=(1,2)$.
		In the second move of Staller, she plays $s_2=(2,1)$ if it is possible. Otherwise, she will play $s_2=(2,3)$. Then Dominator cannot dominate $Z$ within three moves.
	\end{itemize}
	Thus Dominator needs to play at least four moves to win the game which means that $\gmb(Z) \ge 4$.
	
	We conclude that $\gmb (Z) = \gmb' (Z) =4$.
	\qed
	\medskip
	
\begin{figure}[ht]
		\begin{center}
			\begin{tikzpicture}[line cap=round,line join=round,>=triangle 45,x=0.9cm,y=0.9cm]
				\draw [line width=1.2pt] (2.,2.)-- (3.,2.);
				\draw [line width=1.2pt] (3.,2.)-- (4.,2.);
				\draw [line width=1.2pt] (4.,2.)-- (5.,2.);
				\draw [line width=1.2pt] (2.,1.)-- (3.,1.);
				\draw [line width=1.2pt] (3.,1.)-- (4.,1.);
				\draw [line width=1.2pt] (4.,1.)-- (5.,1.);
				\draw [line width=1.2pt] (2.,0.)-- (3.,0.);
				\draw [line width=1.2pt] (3.,0.)-- (4.,0.);
				\draw [line width=1.2pt] (4.,0.)-- (5.,0.);
				\draw [line width=1.2pt] (2.,2.)-- (2.,1.);
				\draw [line width=1.2pt] (2.,1.)-- (2.,0.);
				\draw [line width=1.2pt] (3.,2.)-- (3.,1.);
				\draw [line width=1.2pt] (3.,1.)-- (3.,0.);
				\draw [line width=1.2pt] (4.,2.)-- (4.,1.);
				\draw [line width=1.2pt] (4.,1.)-- (4.,0.);
				\draw [line width=1.2pt] (5.,2.)-- (5.,1.);
				\draw [line width=1.2pt] (5.,1.)-- (5.,0.);
				\draw [line width=1.2pt] (6.,2.)-- (7.,2.);
				\draw [line width=1.2pt] (7.,2.)-- (8.,2.);
				\draw [line width=1.2pt] (8.,2.)-- (9.,2.);
				\draw [line width=1.2pt] (6.,1.)-- (7.,1.);
				\draw [line width=1.2pt] (7.,1.)-- (8.,1.);
				\draw [line width=1.2pt] (8.,1.)-- (9.,1.);
				\draw [line width=1.2pt] (6.,0.)-- (7.,0.);
				\draw [line width=1.2pt] (7.,0.)-- (8.,0.);
				\draw [line width=1.2pt] (8.,0.)-- (9.,0.);
				\draw [line width=1.2pt] (6.,2.)-- (6.,1.);
				\draw [line width=1.2pt] (6.,1.)-- (6.,0.);
				\draw [line width=1.2pt] (7.,2.)-- (7.,1.);
				\draw [line width=1.2pt] (7.,1.)-- (7.,0.);
				\draw [line width=1.2pt] (8.,2.)-- (8.,1.);
				\draw [line width=1.2pt] (8.,1.)-- (8.,0.);
				\draw [line width=1.2pt] (9.,2.)-- (9.,1.);
				\draw [line width=1.2pt] (9.,1.)-- (9.,0.);
				\draw (1.35,0.65) node[anchor=north west] {$d_1$};
				\draw (2.4,2.6) node[anchor=north west] {$s_1$};
				\draw (8.35,0.65) node[anchor=north west] {$d_2$};
				\draw (4.4,0.6) node[anchor=north west] {$s_2$};
				\draw (5.35,0.65) node[anchor=north west] {$d_1$};
				\draw (6.4,2.6) node[anchor=north west] {$s_1$};
				\draw (8.4,2.6) node[anchor=north west] {$s_2$};
				\draw (2.9,-0.1) node[anchor=north west] {$(a)$};
				\draw (6.9,-0.1) node[anchor=north west] {$(b)$};
				\begin{scriptsize}
					\draw [fill=white] (2.,2.) circle (2.5pt);
					\draw [fill=blue] (3.,2.) circle (2.5pt);
					\draw [fill=white] (4.,2.) circle (2.5pt);
					\draw [fill=white] (5.,2.) circle (2.5pt);
					\draw [fill=white] (2.,1.) circle (2.5pt);
					\draw [fill=white] (3.,1.) circle (2.5pt);
					\draw [fill=white] (4.,1.) circle (2.5pt);
					\draw [fill=white] (5.,1.) circle (2.5pt);
					\draw [fill=red] (2.,0.) circle (2.5pt);
					\draw [fill=white] (3.,0.) circle (2.5pt);
					\draw [fill=white] (4.,0.) circle (2.5pt);
					\draw [fill=blue] (5.,0.) circle (2.5pt);
					\draw [fill=white] (6.,2.) circle (2.5pt);
					\draw [fill=blue] (7.,2.) circle (2.5pt);
					\draw [fill=white] (8.,2.) circle (2.5pt);
					\draw [fill=blue] (9.,2.) circle (2.5pt);
					\draw [fill=white] (6.,1.) circle (2.5pt);
					\draw [fill=white] (7.,1.) circle (2.5pt);
					\draw [fill=white] (8.,1.) circle (2.5pt);
					\draw [fill=white] (9.,1.) circle (2.5pt);
					\draw [fill=red] (6.,0.) circle (2.5pt);
					\draw [fill=white] (7.,0.) circle (2.5pt);
					\draw [fill=white] (8.,0.) circle (2.5pt);
					\draw [fill=red] (9.,0.) circle (2.5pt);
				\end{scriptsize}
				
				\draw [line width=1.2pt] (10.,2.)-- (11.,2.);
				\draw [line width=1.2pt] (11.,2.)-- (12.,2.);
				\draw [line width=1.2pt] (12.,2.)-- (13.,2.);
				\draw [line width=1.2pt] (10.,1.)-- (11.,1.);
				\draw [line width=1.2pt] (11.,1.)-- (12.,1.);
				\draw [line width=1.2pt] (12.,1.)-- (13.,1.);
				\draw [line width=1.2pt] (10.,0.)-- (11.,0.);
				\draw [line width=1.2pt] (11.,0.)-- (12.,0.);
				\draw [line width=1.2pt] (12.,0.)-- (13.,0.);
				\draw [line width=1.2pt] (10.,2.)-- (10.,1.);
				\draw [line width=1.2pt] (10.,1.)-- (10.,0.);
				\draw [line width=1.2pt] (11.,2.)-- (11.,1.);
				\draw [line width=1.2pt] (11.,1.)-- (11.,0.);
				\draw [line width=1.2pt] (12.,2.)-- (12.,1.);
				\draw [line width=1.2pt] (12.,1.)-- (12.,0.);
				\draw [line width=1.2pt] (13.,2.)-- (13.,1.);
				\draw [line width=1.2pt] (13.,1.)-- (13.,0.);
				\draw [line width=1.2pt] (14.,2.)-- (15.,2.);
				\draw [line width=1.2pt] (15.,2.)-- (16.,2.);
				\draw [line width=1.2pt] (16.,2.)-- (17.,2.);
				\draw [line width=1.2pt] (14.,1.)-- (15.,1.);
				\draw [line width=1.2pt] (15.,1.)-- (16.,1.);
				\draw [line width=1.2pt] (16.,1.)-- (17.,1.);
				\draw [line width=1.2pt] (14.,0.)-- (15.,0.);
				\draw [line width=1.2pt] (15.,0.)-- (16.,0.);
				\draw [line width=1.2pt] (16.,0.)-- (17.,0.);
				\draw [line width=1.2pt] (14.,2.)-- (14.,1.);
				\draw [line width=1.2pt] (14.,1.)-- (14.,0.);
				\draw [line width=1.2pt] (15.,2.)-- (15.,1.);
				\draw [line width=1.2pt] (15.,1.)-- (15.,0.);
				\draw [line width=1.2pt] (16.,2.)-- (16.,1.);
				\draw [line width=1.2pt] (16.,1.)-- (16.,0.);
				\draw [line width=1.2pt] (17.,2.)-- (17.,1.);
				\draw [line width=1.2pt] (17.,1.)-- (17.,0.);
				\draw (10.35,0.65) node[anchor=north west] {$d_1$};
				\draw (9.4,2.55) node[anchor=north west] {$s_1$};
				\draw (15.35,2.65) node[anchor=north west] {$d_2$};
				\draw (11.4,2.55) node[anchor=north west] {$s_2$};
				\draw (14.35,0.65) node[anchor=north west] {$d_1$};
				\draw (13.4,2.55) node[anchor=north west] {$s_1$};
				\draw (16.4,0.55) node[anchor=north west] {$s_2$};
				\draw (10.9,-0.1) node[anchor=north west] {$(c)$};
				\draw (14.9,-0.1) node[anchor=north west] {$(d)$};
				\begin{scriptsize}
					\draw [fill=blue] (10.,2.) circle (2.5pt);
					\draw [fill=white] (11.,2.) circle (2.5pt);
					\draw [fill=blue] (12.,2.) circle (2.5pt);
					\draw [fill=white] (13.,2.) circle (2.5pt);
					\draw [fill=white] (10.,1.) circle (2.5pt);
					\draw [fill=white] (11.,1.) circle (2.5pt);
					\draw [fill=white] (12.,1.) circle (2.5pt);
					\draw [fill=white] (13.,1.) circle (2.5pt);
					\draw [fill=white] (10.,0.) circle (2.5pt);
					\draw [fill=red] (11.,0.) circle (2.5pt);
					\draw [fill=white] (12.,0.) circle (2.5pt);
					\draw [fill=white] (13.,0.) circle (2.5pt);
					\draw [fill=blue] (14.,2.) circle (2.5pt);
					\draw [fill=white] (15.,2.) circle (2.5pt);
					\draw [fill=red] (16.,2.) circle (2.5pt);
					\draw [fill=white] (17.,2.) circle (2.5pt);
					\draw [fill=white] (14.,1.) circle (2.5pt);
					\draw [fill=white] (15.,1.) circle (2.5pt);
					\draw [fill=white] (16.,1.) circle (2.5pt);
					\draw [fill=white] (17.,1.) circle (2.5pt);
					\draw [fill=white] (14.,0.) circle (2.5pt);
					\draw [fill=red] (15.,0.) circle (2.5pt);
					\draw [fill=white] (16.,0.) circle (2.5pt);
					\draw [fill=blue] (17.,0.) circle (2.5pt);	
				\end{scriptsize}
			\end{tikzpicture}
		\end{center}
		\begin{center}
			\begin{tikzpicture}[line cap=round,line join=round,>=triangle 45,x=0.9cm,y=0.9cm]
				\draw [line width=1.2pt] (2.,2.)-- (3.,2.);
				\draw [line width=1.2pt] (3.,2.)-- (4.,2.);
				\draw [line width=1.2pt] (4.,2.)-- (5.,2.);
				\draw [line width=1.2pt] (2.,1.)-- (3.,1.);
				\draw [line width=1.2pt] (3.,1.)-- (4.,1.);
				\draw [line width=1.2pt] (4.,1.)-- (5.,1.);
				\draw [line width=1.2pt] (2.,0.)-- (3.,0.);
				\draw [line width=1.2pt] (3.,0.)-- (4.,0.);
				\draw [line width=1.2pt] (4.,0.)-- (5.,0.);
				\draw [line width=1.2pt] (2.,2.)-- (2.,1.);
				\draw [line width=1.2pt] (2.,1.)-- (2.,0.);
				\draw [line width=1.2pt] (3.,2.)-- (3.,1.);
				\draw [line width=1.2pt] (3.,1.)-- (3.,0.);
				\draw [line width=1.2pt] (4.,2.)-- (4.,1.);
				\draw [line width=1.2pt] (4.,1.)-- (4.,0.);
				\draw [line width=1.2pt] (5.,2.)-- (5.,1.);
				\draw [line width=1.2pt] (5.,1.)-- (5.,0.);
				\draw [line width=1.2pt] (6.,2.)-- (7.,2.);
				\draw [line width=1.2pt] (7.,2.)-- (8.,2.);
				\draw [line width=1.2pt] (8.,2.)-- (9.,2.);
				\draw [line width=1.2pt] (6.,1.)-- (7.,1.);
				\draw [line width=1.2pt] (7.,1.)-- (8.,1.);
				\draw [line width=1.2pt] (8.,1.)-- (9.,1.);
				\draw [line width=1.2pt] (6.,0.)-- (7.,0.);
				\draw [line width=1.2pt] (7.,0.)-- (8.,0.);
				\draw [line width=1.2pt] (8.,0.)-- (9.,0.);
				\draw [line width=1.2pt] (6.,2.)-- (6.,1.);
				\draw [line width=1.2pt] (6.,1.)-- (6.,0.);
				\draw [line width=1.2pt] (7.,2.)-- (7.,1.);
				\draw [line width=1.2pt] (7.,1.)-- (7.,0.);
				\draw [line width=1.2pt] (8.,2.)-- (8.,1.);
				\draw [line width=1.2pt] (8.,1.)-- (8.,0.);
				\draw [line width=1.2pt] (9.,2.)-- (9.,1.);
				\draw [line width=1.2pt] (9.,1.)-- (9.,0.);
				\draw (1.35,1.65) node[anchor=north west] {$d_1$};
				\draw (3.5,2.6) node[anchor=north west] {$s_2$};
				\draw (7.35,2.65) node[anchor=north west] {$d_2$};
				\draw (3.35,0.6) node[anchor=north west] {$s_1$};
				\draw (5.35,1.65) node[anchor=north west] {$d_1$};
				\draw (7.4,0.6) node[anchor=north west] {$s_1$};
				\draw (8.4,0.6) node[anchor=north west] {$s_2$};
				\draw (2.9,-0.1) node[anchor=north west] {$(e)$};
				\draw (6.9,-0.1) node[anchor=north west] {$(f)$};
				\begin{scriptsize}
					\draw [fill=white] (2.,2.) circle (2.5pt);
					\draw [fill=white] (3.,2.) circle (2.5pt);
					\draw [fill=blue] (4.,2.) circle (2.5pt);
					\draw [fill=white] (5.,2.) circle (2.5pt);
					\draw [fill=red] (2.,1.) circle (2.5pt);
					\draw [fill=white] (3.,1.) circle (2.5pt);
					\draw [fill=white] (4.,1.) circle (2.5pt);
					\draw [fill=white] (5.,1.) circle (2.5pt);
					\draw [fill=white] (2.,0.) circle (2.5pt);
					\draw [fill=white] (3.,0.) circle (2.5pt);
					\draw [fill=blue] (4.,0.) circle (2.5pt);
					\draw [fill=white] (5.,0.) circle (2.5pt);
					\draw [fill=white] (6.,2.) circle (2.5pt);
					\draw [fill=white] (7.,2.) circle (2.5pt);
					\draw [fill=red] (8.,2.) circle (2.5pt);
					\draw [fill=white] (9.,2.) circle (2.5pt);
					\draw [fill=red] (6.,1.) circle (2.5pt);
					\draw [fill=white] (7.,1.) circle (2.5pt);
					\draw [fill=white] (8.,1.) circle (2.5pt);
					\draw [fill=white] (9.,1.) circle (2.5pt);
					\draw [fill=white] (6.,0.) circle (2.5pt);
					\draw [fill=white] (7.,0.) circle (2.5pt);
					\draw [fill=blue] (8.,0.) circle (2.5pt);
					\draw [fill=blue] (9.,0.) circle (2.5pt);
				\end{scriptsize}
				
				\draw [line width=1.2pt] (10.,2.)-- (11.,2.);
				\draw [line width=1.2pt] (11.,2.)-- (12.,2.);
				\draw [line width=1.2pt] (12.,2.)-- (13.,2.);
				\draw [line width=1.2pt] (10.,1.)-- (11.,1.);
				\draw [line width=1.2pt] (11.,1.)-- (12.,1.);
				\draw [line width=1.2pt] (12.,1.)-- (13.,1.);
				\draw [line width=1.2pt] (10.,0.)-- (11.,0.);
				\draw [line width=1.2pt] (11.,0.)-- (12.,0.);
				\draw [line width=1.2pt] (12.,0.)-- (13.,0.);
				\draw [line width=1.2pt] (10.,2.)-- (10.,1.);
				\draw [line width=1.2pt] (10.,1.)-- (10.,0.);
				\draw [line width=1.2pt] (11.,2.)-- (11.,1.);
				\draw [line width=1.2pt] (11.,1.)-- (11.,0.);
				\draw [line width=1.2pt] (12.,2.)-- (12.,1.);
				\draw [line width=1.2pt] (12.,1.)-- (12.,0.);
				\draw [line width=1.2pt] (13.,2.)-- (13.,1.);
				\draw [line width=1.2pt] (13.,1.)-- (13.,0.);
				\draw [line width=1.2pt] (14.,2.)-- (15.,2.);
				\draw [line width=1.2pt] (15.,2.)-- (16.,2.);
				\draw [line width=1.2pt] (16.,2.)-- (17.,2.);
				\draw [line width=1.2pt] (14.,1.)-- (15.,1.);
				\draw [line width=1.2pt] (15.,1.)-- (16.,1.);
				\draw [line width=1.2pt] (16.,1.)-- (17.,1.);
				\draw [line width=1.2pt] (14.,0.)-- (15.,0.);
				\draw [line width=1.2pt] (15.,0.)-- (16.,0.);
				\draw [line width=1.2pt] (16.,0.)-- (17.,0.);
				\draw [line width=1.2pt] (14.,2.)-- (14.,1.);
				\draw [line width=1.2pt] (14.,1.)-- (14.,0.);
				\draw [line width=1.2pt] (15.,2.)-- (15.,1.);
				\draw [line width=1.2pt] (15.,1.)-- (15.,0.);
				\draw [line width=1.2pt] (16.,2.)-- (16.,1.);
				\draw [line width=1.2pt] (16.,1.)-- (16.,0.);
				\draw [line width=1.2pt] (17.,2.)-- (17.,1.);
				\draw [line width=1.2pt] (17.,1.)-- (17.,0.);
				\draw (10.35,1.65) node[anchor=north west] {$d_1$};
				\draw (9.4,1.6) node[anchor=north west] {$s_1$};
				\draw (16.35,0.65) node[anchor=north west] {$d_2$};
				\draw (12.4,0.6) node[anchor=north west] {$s_2$};
				\draw (14.4,1.65) node[anchor=north west] {$d_1$};
				\draw (13.4,1.6) node[anchor=north west] {$s_1$};
				\draw (16.4,2.6) node[anchor=north west] {$s_2$};
				\draw (10.9,-0.1) node[anchor=north west] {$(g)$};
				\draw (14.9,-0.1) node[anchor=north west] {$(h)$};
				\begin{scriptsize}
					\draw [fill=white] (10.,2.) circle (2.5pt);
					\draw [fill=white] (11.,2.) circle (2.5pt);
					\draw [fill=white] (12.,2.) circle (2.5pt);
					\draw [fill=white] (13.,2.) circle (2.5pt);
					\draw [fill=blue] (10.,1.) circle (2.5pt);
					\draw [fill=red] (11.,1.) circle (2.5pt);
					\draw [fill=white] (12.,1.) circle (2.5pt);
					\draw [fill=white] (13.,1.) circle (2.5pt);
					\draw [fill=white] (10.,0.) circle (2.5pt);
					\draw [fill=white] (11.,0.) circle (2.5pt);
					\draw [fill=white] (12.,0.) circle (2.5pt);
					\draw [fill=blue] (13.,0.) circle (2.5pt);
					\draw [fill=white] (14.,2.) circle (2.5pt);
					\draw [fill=white] (15.,2.) circle (2.5pt);
					\draw [fill=white] (16.,2.) circle (2.5pt);
					\draw [fill=blue] (17.,2.) circle (2.5pt);
					\draw [fill=blue] (14.,1.) circle (2.5pt);
					\draw [fill=red] (15.,1.) circle (2.5pt);
					\draw [fill=white] (16.,1.) circle (2.5pt);
					\draw [fill=white] (17.,1.) circle (2.5pt);
					\draw [fill=white] (14.,0.) circle (2.5pt);
					\draw [fill=white] (15.,0.) circle (2.5pt);
					\draw [fill=white] (16.,0.) circle (2.5pt);
					\draw [fill=red] (17.,0.) circle (2.5pt);	
				\end{scriptsize}
			\end{tikzpicture}
			\caption{An illustration for Staller's strategy in the proof of $\gmb(P_3 \cp P_4) \ge 5$. Blue and red vertices denote the moves of Staller and Dominator, respectively.}
			\label{fig:Staller for P3P4}
		\end{center}
	\end{figure}
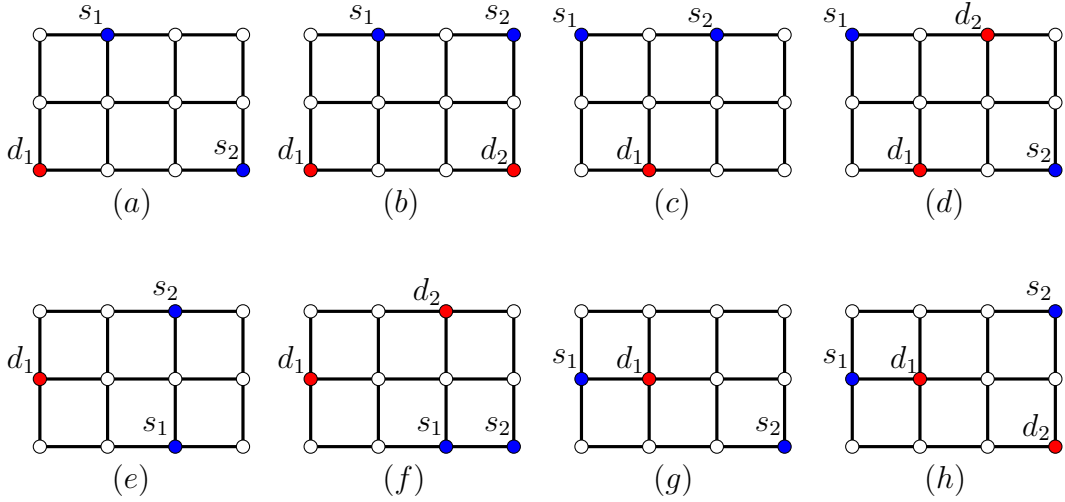

	\begin{proposition}
		$\gmb(P_3 \cp P_4) =5$ and $\gmb'(P_3 \cp P_4) =6$.
	\end{proposition}
	\proof 
	Let $Z = P_3 \cp P_4$ and we first consider the D-game.
	We will prove that $\gmb(Z) \ge 5$ by providing a strategy for Staller which ensures that Dominator cannot form a dominating set within four moves in the D-game. Consider the first move $d_1$ of Dominator.
	\begin{itemize}
		\item  If $d_1 = (1,1)$, then Staller plays $s_1=(3,2)$. In her second move, she plays $s_2 = (1,4)$ if it is possible. Otherwise she plays $s_2= (3,4)$. See Figure \ref{fig:Staller for P3P4}\,(a),\,(b).
		\item  If $d_1 = (1,2)$, then Staller replies by playing $s_1=(3,1)$. After that, she plays $s_2=(3,3)$ if it is possible, otherwise she plays $s_2=(1,4)$.  See Figure \ref{fig:Staller for P3P4}\,(c),\,(d).
		\item  If $d_1 = (2,1)$, then Staller responds at $s_1=(1,3)$. Then she plays $s_2=(3,3)$ if it is possible, otherwise she plays $s_2=(1,4)$. See Figure \ref{fig:Staller for P3P4}\,(e),\,(f).
		\item  If $d_1 = (2,2)$, then Staller replies by playing $s_1=(2,1)$. In her second move, she plays $s_2 \in \{(1,4), (3,4)\}$. See Figure \ref{fig:Staller for P3P4}\,(g),\,(h).
	\end{itemize}
	From the above strategies, one can see that Dominator cannot form a dominating set within four moves. Therefore $\gmb(Z)  \ge 5.$
	
	Next, let $Z'$ be a graph obtained from $Z$ by deleting $(i,1)(i,2)$ where $i \in [3]$. By Lemma~\ref{lem:o_edgedelete}\, and Theorem~\ref{thm:pmpn}, it follows that $\gmb(Z) \le \gmb(Z')$. We now claim that $\gmb(Z') \le 5$.
	Observe that $Z'$ is a disjoint union of $P_3$ and $P_3 \cp P_3$.  Assume that Dominator starts the game by playing $d_1=(2,1)$. Then he dominated the component $P_3$ and both players will continue the game in $P_3\cp P_3$. By Proposition~\ref{prop:p3p3}, Dominator needs four more move to win the game in $Z'$. It concludes that $\gmb(Z) \le \gmb(Z') \le 5$.

	\medskip
	Now consider the S-game. Since Dominator is the winner in the S-game, $\gmb'(Z) \le 6$ by (\ref{bound2}).
	It remains to show that $\gmb'(Z) \ge 6$. We will provide a strategy for Staller to ensure that Dominator needs to play at least six moves. Staller starts the game with $s'_1 =(2,1)$.
	
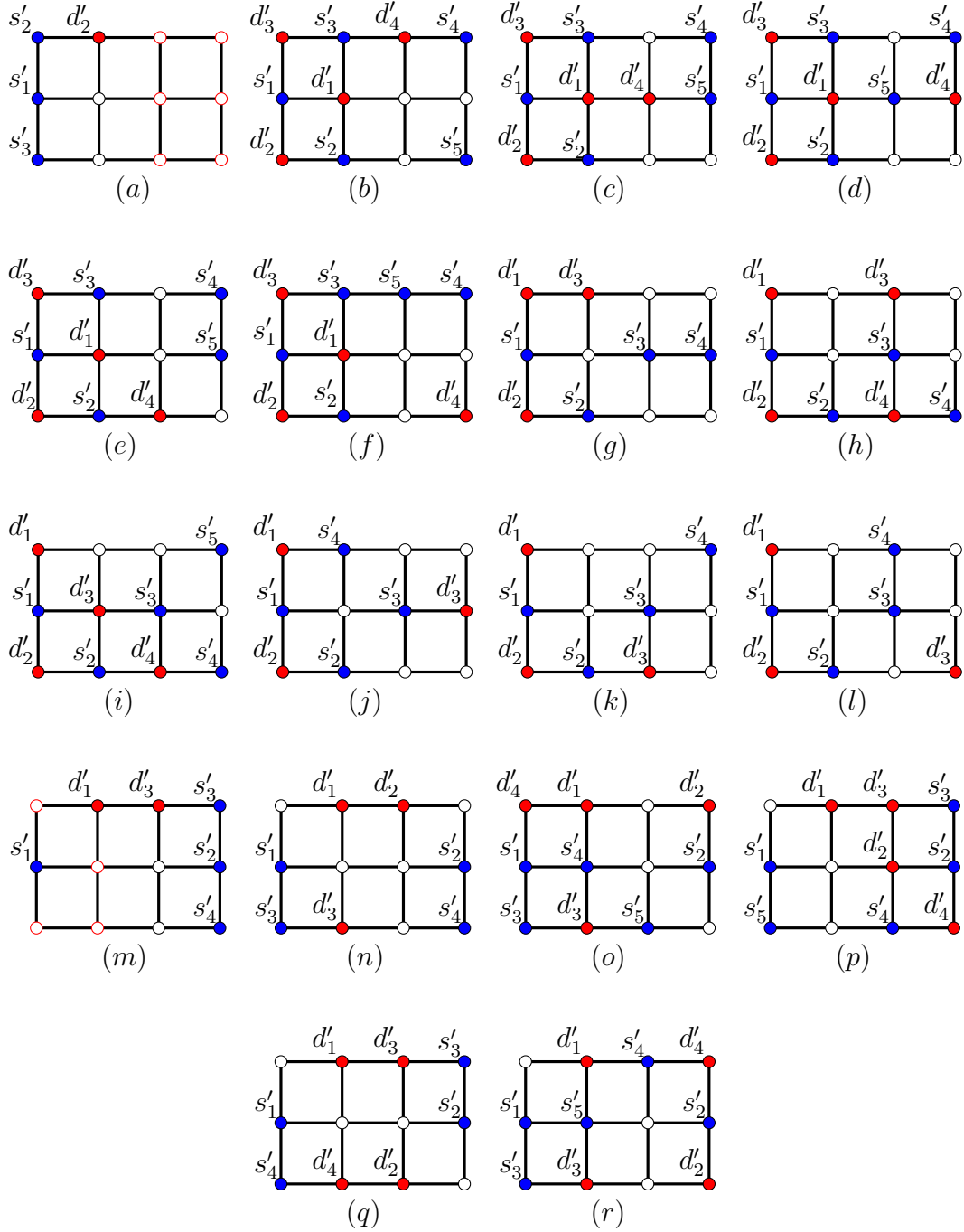
\begin{figure}[H]
		
		\begin{center}
			\begin{tikzpicture}[line cap=round,line join=round,>=triangle 45,x=0.9cm,y=0.9cm]
				\draw [line width=1.2pt] (2.,2.)-- (3.,2.);
				\draw [line width=1.2pt] (3.,2.)-- (4.,2.);
				\draw [line width=1.2pt] (4.,2.)-- (5.,2.);
				\draw [line width=1.2pt] (2.,1.)-- (3.,1.);
				\draw [line width=1.2pt] (3.,1.)-- (4.,1.);
				\draw [line width=1.2pt] (4.,1.)-- (5.,1.);
				\draw [line width=1.2pt] (2.,0.)-- (3.,0.);
				\draw [line width=1.2pt] (3.,0.)-- (4.,0.);
				\draw [line width=1.2pt] (4.,0.)-- (5.,0.);
				\draw [line width=1.2pt] (2.,2.)-- (2.,1.);
				\draw [line width=1.2pt] (2.,1.)-- (2.,0.);
				\draw [line width=1.2pt] (3.,2.)-- (3.,1.);
				\draw [line width=1.2pt] (3.,1.)-- (3.,0.);
				\draw [line width=1.2pt] (4.,2.)-- (4.,1.);
				\draw [line width=1.2pt] (4.,1.)-- (4.,0.);
				\draw [line width=1.2pt] (5.,2.)-- (5.,1.);
				\draw [line width=1.2pt] (5.,1.)-- (5.,0.);
				\draw [line width=1.2pt] (6.,2.)-- (7.,2.);
				\draw [line width=1.2pt] (7.,2.)-- (8.,2.);
				\draw [line width=1.2pt] (8.,2.)-- (9.,2.);
				\draw [line width=1.2pt] (6.,1.)-- (7.,1.);
				\draw [line width=1.2pt] (7.,1.)-- (8.,1.);
				\draw [line width=1.2pt] (8.,1.)-- (9.,1.);
				\draw [line width=1.2pt] (6.,0.)-- (7.,0.);
				\draw [line width=1.2pt] (7.,0.)-- (8.,0.);
				\draw [line width=1.2pt] (8.,0.)-- (9.,0.);
				\draw [line width=1.2pt] (6.,2.)-- (6.,1.);
				\draw [line width=1.2pt] (6.,1.)-- (6.,0.);
				\draw [line width=1.2pt] (7.,2.)-- (7.,1.);
				\draw [line width=1.2pt] (7.,1.)-- (7.,0.);
				\draw [line width=1.2pt] (8.,2.)-- (8.,1.);
				\draw [line width=1.2pt] (8.,1.)-- (8.,0.);
				\draw [line width=1.2pt] (9.,2.)-- (9.,1.);
				\draw [line width=1.2pt] (9.,1.)-- (9.,0.);
				\draw (6.3,1.7) node[anchor=north west] {$d'_1$};
				\draw (6.35,0.7) node[anchor=north west] {$s'_2$};
				\draw (5.3,0.7) node[anchor=north west] {$d'_2$};
				\draw (1.35,1.7) node[anchor=north west] {$s'_1$};
				\draw (6.35,2.7) node[anchor=north west] {$s'_3$};
				\draw (5.3,2.7) node[anchor=north west] {$d'_3$};
				\draw (2.3,2.7) node[anchor=north west] {$d'_2$};
				\draw (1.35,2.7) node[anchor=north west] {$s'_2$};
				\draw (1.35,0.7) node[anchor=north west] {$s'_3$};
				\draw (5.35,1.7) node[anchor=north west] {$s'_1$};
				\draw (8.35,2.7) node[anchor=north west] {$s'_4$};
				\draw (3.1,-0.1) node[anchor=north west] { $(a)$};
				\draw (6.9,-0.1) node[anchor=north west] { $(b)$};
				\begin{scriptsize}
					\draw [fill=blue] (2.,2.) circle (2.5pt);
					\draw [fill=red] (3.,2.) circle (2.5pt);
					\draw [color=red, fill=white] (4.,2.) circle (2.5pt);
					\draw [color=red, fill=white] (5.,2.) circle (2.5pt);
					\draw [fill=blue] (2.,1.) circle (2.5pt);
					\draw [fill=white] (3.,1.) circle (2.5pt);
					\draw [color=red, fill=white] (4.,1.) circle (2.5pt);
					\draw [color=red, fill=white] (5.,1.) circle (2.5pt);
					\draw [fill=blue] (2.,0.) circle (2.5pt);
					\draw [fill=white] (3.,0.) circle (2.5pt);
					\draw [color=red, fill=white] (4.,0.) circle (2.5pt);
					\draw [color=red, fill=white] (5.,0.) circle (2.5pt);
					\draw [fill=red] (6.,2.) circle (2.5pt);
					\draw [fill=blue] (7.,2.) circle (2.5pt);
					\draw [fill=red] (8.,2.) circle (2.5pt);
					\draw [fill=blue] (9.,2.) circle (2.5pt);
					\draw [fill=blue] (6.,1.) circle (2.5pt);
					\draw [fill=red] (7.,1.) circle (2.5pt);
					\draw [fill=white] (8.,1.) circle (2.5pt);
					\draw [fill=white] (9.,1.) circle (2.5pt);
					\draw [fill=red] (6.,0.) circle (2.5pt);
					\draw [fill=blue] (7.,0.) circle (2.5pt);
					\draw [fill=white] (8.,0.) circle (2.5pt);
					\draw [fill=blue] (9.,0.) circle (2.5pt);
				\end{scriptsize}
				\draw [line width=1.2pt] (10.,2.)-- (11.,2.);
				\draw [line width=1.2pt] (11.,2.)-- (12.,2.);
				\draw [line width=1.2pt] (12.,2.)-- (13.,2.);
				\draw [line width=1.2pt] (10.,1.)-- (11.,1.);
				\draw [line width=1.2pt] (11.,1.)-- (12.,1.);
				\draw [line width=1.2pt] (12.,1.)-- (13.,1.);
				\draw [line width=1.2pt] (10.,0.)-- (11.,0.);
				\draw [line width=1.2pt] (11.,0.)-- (12.,0.);
				\draw [line width=1.2pt] (12.,0.)-- (13.,0.);
				\draw [line width=1.2pt] (10.,2.)-- (10.,1.);
				\draw [line width=1.2pt] (10.,1.)-- (10.,0.);
				\draw [line width=1.2pt] (11.,2.)-- (11.,1.);
				\draw [line width=1.2pt] (11.,1.)-- (11.,0.);
				\draw [line width=1.2pt] (12.,2.)-- (12.,1.);
				\draw [line width=1.2pt] (12.,1.)-- (12.,0.);
				\draw [line width=1.2pt] (13.,2.)-- (13.,1.);
				\draw [line width=1.2pt] (13.,1.)-- (13.,0.);
				\draw [line width=1.2pt] (14.,2.)-- (15.,2.);
				\draw [line width=1.2pt] (15.,2.)-- (16.,2.);
				\draw [line width=1.2pt] (16.,2.)-- (17.,2.);
				\draw [line width=1.2pt] (14.,1.)-- (15.,1.);
				\draw [line width=1.2pt] (15.,1.)-- (16.,1.);
				\draw [line width=1.2pt] (16.,1.)-- (17.,1.);
				\draw [line width=1.2pt] (14.,0.)-- (15.,0.);
				\draw [line width=1.2pt] (15.,0.)-- (16.,0.);
				\draw [line width=1.2pt] (16.,0.)-- (17.,0.);
				\draw [line width=1.2pt] (14.,2.)-- (14.,1.);
				\draw [line width=1.2pt] (14.,1.)-- (14.,0.);
				\draw [line width=1.2pt] (15.,2.)-- (15.,1.);
				\draw [line width=1.2pt] (15.,1.)-- (15.,0.);
				\draw [line width=1.2pt] (16.,2.)-- (16.,1.);
				\draw [line width=1.2pt] (16.,1.)-- (16.,0.);
				\draw [line width=1.2pt] (17.,2.)-- (17.,1.);
				\draw [line width=1.2pt] (17.,1.)-- (17.,0.);
				\draw (10.35,1.75) node[anchor=north west] {$d'_1$};
				\draw (10.4,0.65) node[anchor=north west] {$s'_2$};
				\draw (9.35,0.75) node[anchor=north west] {$d'_2$};
				\draw (9.4,1.7) node[anchor=north west] {$s'_1$};
				\draw (10.4,2.7) node[anchor=north west] {$s'_3$};
				\draw (9.4,2.75) node[anchor=north west] {$d'_3$};
				\draw (12.4,2.7) node[anchor=north west] {$s'_4$};
				\draw (7.35,2.75) node[anchor=north west] {$d'_4$};
				\draw (8.4,0.7) node[anchor=north west] {$s'_5$};
				\draw (11.35,1.75) node[anchor=north west] {$d'_4$};
				\draw (12.4,1.7) node[anchor=north west] {$s'_5$};
				\draw (14.35,1.75) node[anchor=north west] {$d'_1$};
				\draw (14.4,0.7) node[anchor=north west] {$s'_2$};
				\draw (13.35,0.75) node[anchor=north west] {$d'_2$};
				\draw (13.4,1.7) node[anchor=north west] {$s'_1$};
				\draw (14.4,2.7) node[anchor=north west] {$s'_3$};
				\draw (13.35,2.75) node[anchor=north west] {$d'_3$};
				\draw (16.4,2.7) node[anchor=north west] {$s'_4$};
				\draw (16.35,1.75) node[anchor=north west] {$d'_4$};
				\draw (15.4,1.7) node[anchor=north west] {$s'_5$};
				\draw (10.9,-0.1) node[anchor=north west] {$(c)$};
				\draw (14.9,-0.1) node[anchor=north west] {$(d)$};
				\begin{scriptsize}
					\draw [fill=red] (10.,2.) circle (2.5pt);
					\draw [fill=blue] (11.,2.) circle (2.5pt);
					\draw [fill=white] (12.,2.) circle (2.5pt);
					\draw [fill=blue] (13.,2.) circle (2.5pt);
					\draw [fill=blue] (10.,1.) circle (2.5pt);
					\draw [fill=red] (11.,1.) circle (2.5pt);
					\draw [fill=red] (12.,1.) circle (2.5pt);
					\draw [fill=blue] (13.,1.) circle (2.5pt);
					\draw [fill=red] (10.,0.) circle (2.5pt);
					\draw [fill=blue] (11.,0.) circle (2.5pt);
					\draw [fill=white] (12.,0.) circle (2.5pt);
					\draw [fill=white] (13.,0.) circle (2.5pt);
					\draw [fill=red] (14.,2.) circle (2.5pt);
					\draw [fill=blue] (15.,2.) circle (2.5pt);
					\draw [fill=white] (16.,2.) circle (2.5pt);
					\draw [fill=blue] (17.,2.) circle (2.5pt);
					\draw [fill=blue] (14.,1.) circle (2.5pt);
					\draw [fill=red] (15.,1.) circle (2.5pt);
					\draw [fill=blue] (16.,1.) circle (2.5pt);
					\draw [fill=red] (17.,1.) circle (2.5pt);
					\draw [fill=red] (14.,0.) circle (2.5pt);
					\draw [fill=blue] (15.,0.) circle (2.5pt);
					\draw [fill=white] (16.,0.) circle (2.5pt);
					\draw [fill=white] (17.,0.) circle (2.5pt);	
				\end{scriptsize}
			\end{tikzpicture}
		\end{center}	
		\begin{center}
			\begin{tikzpicture}[line cap=round,line join=round,>=triangle 45,x=0.9cm,y=0.9cm]
				\draw [line width=1.2pt] (2.,2.)-- (3.,2.);
				\draw [line width=1.2pt] (3.,2.)-- (4.,2.);
				\draw [line width=1.2pt] (4.,2.)-- (5.,2.);
				\draw [line width=1.2pt] (2.,1.)-- (3.,1.);
				\draw [line width=1.2pt] (3.,1.)-- (4.,1.);
				\draw [line width=1.2pt] (4.,1.)-- (5.,1.);
				\draw [line width=1.2pt] (2.,0.)-- (3.,0.);
				\draw [line width=1.2pt] (3.,0.)-- (4.,0.);
				\draw [line width=1.2pt] (4.,0.)-- (5.,0.);
				\draw [line width=1.2pt] (2.,2.)-- (2.,1.);
				\draw [line width=1.2pt] (2.,1.)-- (2.,0.);
				\draw [line width=1.2pt] (3.,2.)-- (3.,1.);
				\draw [line width=1.2pt] (3.,1.)-- (3.,0.);
				\draw [line width=1.2pt] (4.,2.)-- (4.,1.);
				\draw [line width=1.2pt] (4.,1.)-- (4.,0.);
				\draw [line width=1.2pt] (5.,2.)-- (5.,1.);
				\draw [line width=1.2pt] (5.,1.)-- (5.,0.);
				\draw [line width=1.2pt] (6.,2.)-- (7.,2.);
				\draw [line width=1.2pt] (7.,2.)-- (8.,2.);
				\draw [line width=1.2pt] (8.,2.)-- (9.,2.);
				\draw [line width=1.2pt] (6.,1.)-- (7.,1.);
				\draw [line width=1.2pt] (7.,1.)-- (8.,1.);
				\draw [line width=1.2pt] (8.,1.)-- (9.,1.);
				\draw [line width=1.2pt] (6.,0.)-- (7.,0.);
				\draw [line width=1.2pt] (7.,0.)-- (8.,0.);
				\draw [line width=1.2pt] (8.,0.)-- (9.,0.);
				\draw [line width=1.2pt] (6.,2.)-- (6.,1.);
				\draw [line width=1.2pt] (6.,1.)-- (6.,0.);
				\draw [line width=1.2pt] (7.,2.)-- (7.,1.);
				\draw [line width=1.2pt] (7.,1.)-- (7.,0.);
				\draw [line width=1.2pt] (8.,2.)-- (8.,1.);
				\draw [line width=1.2pt] (8.,1.)-- (8.,0.);
				\draw [line width=1.2pt] (9.,2.)-- (9.,1.);
				\draw [line width=1.2pt] (9.,1.)-- (9.,0.);
				
				\draw (6.35,1.75) node[anchor=north west] {$d'_1$};
				\draw (2.35,1.75) node[anchor=north west] {$d'_1$};
				\draw (6.35,0.77) node[anchor=north west] {$s'_2$};
				\draw (2.4,0.7) node[anchor=north west] {$s'_2$};
				\draw (5.35,0.72) node[anchor=north west] {$d'_2$};
				\draw (1.4,1.7) node[anchor=north west] {$s'_1$};
				\draw (6.4,2.7) node[anchor=north west] {$s'_3$};
				\draw (5.35,2.75) node[anchor=north west] {$d'_3$};
				\draw (2.4,2.7) node[anchor=north west] {$s'_3$};
				\draw (1.35,2.75) node[anchor=north west] {$d'_3$};
				\draw (1.4,0.72) node[anchor=north west] {$d'_2$};
				\draw (5.35,1.77) node[anchor=north west] {$s'_1$};
				\draw (4.4,2.7) node[anchor=north west] {$s'_4$};
				\draw (8.4,2.7) node[anchor=north west] {$s'_4$};
				\draw (3.35,0.75) node[anchor=north west] {$d'_4$};
				\draw (4.4,1.7) node[anchor=north west] {$s'_5$};
				\draw (8.35,0.75) node[anchor=north west] {$d'_4$};
				\draw (7.4,2.7) node[anchor=north west] {$s'_5$};
				\draw (2.9,-0.1) node[anchor=north west] { $(e)$};
				\draw (6.9,-0.1) node[anchor=north west] { $(f)$};
				\begin{scriptsize}
					\draw [fill=red] (2.,2.) circle (2.5pt);
					\draw [fill=blue] (3.,2.) circle (2.5pt);
					\draw [fill=white] (4.,2.) circle (2.5pt);
					\draw [fill=blue] (5.,2.) circle (2.5pt);
					\draw [fill=blue] (2.,1.) circle (2.5pt);
					\draw [fill=red] (3.,1.) circle (2.5pt);
					\draw [fill=white] (4.,1.) circle (2.5pt);
					\draw [fill=blue] (5.,1.) circle (2.5pt);
					\draw [fill=red] (2.,0.) circle (2.5pt);
					\draw [fill=blue] (3.,0.) circle (2.5pt);
					\draw [fill=red] (4.,0.) circle (2.5pt);
					\draw [fill=white] (5.,0.) circle (2.5pt);
					\draw [fill=red] (6.,2.) circle (2.5pt);
					\draw [fill=blue] (7.,2.) circle (2.5pt);
					\draw [fill=blue] (8.,2.) circle (2.5pt);
					\draw [fill=blue] (9.,2.) circle (2.5pt);
					\draw [fill=blue] (6.,1.) circle (2.5pt);
					\draw [fill=red] (7.,1.) circle (2.5pt);
					\draw [fill=white] (8.,1.) circle (2.5pt);
					\draw [fill=white] (9.,1.) circle (2.5pt);
					\draw [fill=red] (6.,0.) circle (2.5pt);
					\draw [fill=blue] (7.,0.) circle (2.5pt);
					\draw [fill=white] (8.,0.) circle (2.5pt);
					\draw [fill=red] (9.,0.) circle (2.5pt);
				\end{scriptsize}
				
				\draw [line width=1.2pt] (10.,2.)-- (11.,2.);
				\draw [line width=1.2pt] (11.,2.)-- (12.,2.);
				\draw [line width=1.2pt] (12.,2.)-- (13.,2.);
				\draw [line width=1.2pt] (10.,1.)-- (11.,1.);
				\draw [line width=1.2pt] (11.,1.)-- (12.,1.);
				\draw [line width=1.2pt] (12.,1.)-- (13.,1.);
				\draw [line width=1.2pt] (10.,0.)-- (11.,0.);
				\draw [line width=1.2pt] (11.,0.)-- (12.,0.);
				\draw [line width=1.2pt] (12.,0.)-- (13.,0.);
				\draw [line width=1.2pt] (10.,2.)-- (10.,1.);
				\draw [line width=1.2pt] (10.,1.)-- (10.,0.);
				\draw [line width=1.2pt] (11.,2.)-- (11.,1.);
				\draw [line width=1.2pt] (11.,1.)-- (11.,0.);
				\draw [line width=1.2pt] (12.,2.)-- (12.,1.);
				\draw [line width=1.2pt] (12.,1.)-- (12.,0.);
				\draw [line width=1.2pt] (13.,2.)-- (13.,1.);
				\draw [line width=1.2pt] (13.,1.)-- (13.,0.);
				\draw [line width=1.2pt] (14.,2.)-- (15.,2.);
				\draw [line width=1.2pt] (15.,2.)-- (16.,2.);
				\draw [line width=1.2pt] (16.,2.)-- (17.,2.);
				\draw [line width=1.2pt] (14.,1.)-- (15.,1.);
				\draw [line width=1.2pt] (15.,1.)-- (16.,1.);
				\draw [line width=1.2pt] (16.,1.)-- (17.,1.);
				\draw [line width=1.2pt] (14.,0.)-- (15.,0.);
				\draw [line width=1.2pt] (15.,0.)-- (16.,0.);
				\draw [line width=1.2pt] (16.,0.)-- (17.,0.);
				\draw [line width=1.2pt] (14.,2.)-- (14.,1.);
				\draw [line width=1.2pt] (14.,1.)-- (14.,0.);
				\draw [line width=1.2pt] (15.,2.)-- (15.,1.);
				\draw [line width=1.2pt] (15.,1.)-- (15.,0.);
				\draw [line width=1.2pt] (16.,2.)-- (16.,1.);
				\draw [line width=1.2pt] (16.,1.)-- (16.,0.);
				\draw [line width=1.2pt] (17.,2.)-- (17.,1.);
				\draw [line width=1.2pt] (17.,1.)-- (17.,0.);
				\draw (9.4,1.7) node[anchor=north west] {$s'_1$};
				\draw (9.35,2.75) node[anchor=north west] {$d'_1$};
				\draw (10.4,0.7) node[anchor=north west] {$s'_2$};
				\draw (9.35,0.75) node[anchor=north west] {$d'_2$};
				\draw (11.4,1.7) node[anchor=north west] {$s'_3$};
				\draw (12.4,1.7) node[anchor=north west] {$s'_4$};
				\draw (10.35,2.75) node[anchor=north west] {$d'_3$};
				
				\draw (13.4,1.7) node[anchor=north west] {$s'_1$};
				\draw (13.35,2.75) node[anchor=north west] {$d'_1$};
				\draw (14.4,0.7) node[anchor=north west] {$s'_2$};
				\draw (13.35,0.75) node[anchor=north west] {$d'_2$};
				\draw (15.4,1.7) node[anchor=north west] {$s'_3$};
				\draw (16.4,0.7) node[anchor=north west] {$s'_4$};
				\draw (15.35,2.75) node[anchor=north west] {$d'_3$};
				\draw (15.35,0.75) node[anchor=north west] {$d'_4$};
				\draw (10.9,-0.1) node[anchor=north west] {$(g)$};
				\draw (14.9,-0.1) node[anchor=north west] {$(h)$};
				
				\begin{scriptsize}
					\draw [fill=red] (10.,2.) circle (2.5pt);
					\draw [fill=red] (11.,2.) circle (2.5pt);
					\draw [fill=white] (12.,2.) circle (2.5pt);
					\draw [fill=white] (13.,2.) circle (2.5pt);
					\draw [fill=blue] (10.,1.) circle (2.5pt);
					\draw [fill=white] (11.,1.) circle (2.5pt);
					\draw [fill=blue] (12.,1.) circle (2.5pt);
					\draw [fill=blue] (13.,1.) circle (2.5pt);
					\draw [fill=red] (10.,0.) circle (2.5pt);
					\draw [fill=blue] (11.,0.) circle (2.5pt);
					\draw [fill=white] (12.,0.) circle (2.5pt);
					\draw [fill=white] (13.,0.) circle (2.5pt);
					\draw [fill=red] (14.,2.) circle (2.5pt);
					\draw [fill=white] (15.,2.) circle (2.5pt);
					\draw [fill=red] (16.,2.) circle (2.5pt);
					\draw [fill=white] (17.,2.) circle (2.5pt);
					\draw [fill=blue] (14.,1.) circle (2.5pt);
					\draw [fill=white] (15.,1.) circle (2.5pt);
					\draw [fill=blue] (16.,1.) circle (2.5pt);
					\draw [fill=white] (17.,1.) circle (2.5pt);
					\draw [fill=red] (14.,0.) circle (2.5pt);
					\draw [fill=blue] (15.,0.) circle (2.5pt);
					\draw [fill=red] (16.,0.) circle (2.5pt);
					\draw [fill=blue] (17.,0.) circle (2.5pt);	
				\end{scriptsize}
			\end{tikzpicture}
		\end{center}
		\begin{center}
			\begin{tikzpicture}[line cap=round,line join=round,>=triangle 45,x=0.9cm,y=0.9cm]
				\draw [line width=1.2pt] (2.,2.)-- (3.,2.);
				\draw [line width=1.2pt] (3.,2.)-- (4.,2.);
				\draw [line width=1.2pt] (4.,2.)-- (5.,2.);
				\draw [line width=1.2pt] (2.,1.)-- (3.,1.);
				\draw [line width=1.2pt] (3.,1.)-- (4.,1.);
				\draw [line width=1.2pt] (4.,1.)-- (5.,1.);
				\draw [line width=1.2pt] (2.,0.)-- (3.,0.);
				\draw [line width=1.2pt] (3.,0.)-- (4.,0.);
				\draw [line width=1.2pt] (4.,0.)-- (5.,0.);
				\draw [line width=1.2pt] (2.,2.)-- (2.,1.);
				\draw [line width=1.2pt] (2.,1.)-- (2.,0.);
				\draw [line width=1.2pt] (3.,2.)-- (3.,1.);
				\draw [line width=1.2pt] (3.,1.)-- (3.,0.);
				\draw [line width=1.2pt] (4.,2.)-- (4.,1.);
				\draw [line width=1.2pt] (4.,1.)-- (4.,0.);
				\draw [line width=1.2pt] (5.,2.)-- (5.,1.);
				\draw [line width=1.2pt] (5.,1.)-- (5.,0.);
				\draw [line width=1.2pt] (6.,2.)-- (7.,2.);
				\draw [line width=1.2pt] (7.,2.)-- (8.,2.);
				\draw [line width=1.2pt] (8.,2.)-- (9.,2.);
				\draw [line width=1.2pt] (6.,1.)-- (7.,1.);
				\draw [line width=1.2pt] (7.,1.)-- (8.,1.);
				\draw [line width=1.2pt] (8.,1.)-- (9.,1.);
				\draw [line width=1.2pt] (6.,0.)-- (7.,0.);
				\draw [line width=1.2pt] (7.,0.)-- (8.,0.);
				\draw [line width=1.2pt] (8.,0.)-- (9.,0.);
				\draw [line width=1.2pt] (6.,2.)-- (6.,1.);
				\draw [line width=1.2pt] (6.,1.)-- (6.,0.);
				\draw [line width=1.2pt] (7.,2.)-- (7.,1.);
				\draw [line width=1.2pt] (7.,1.)-- (7.,0.);
				\draw [line width=1.2pt] (8.,2.)-- (8.,1.);
				\draw [line width=1.2pt] (8.,1.)-- (8.,0.);
				\draw [line width=1.2pt] (9.,2.)-- (9.,1.);
				\draw [line width=1.2pt] (9.,1.)-- (9.,0.);
				\draw (1.4,1.7) node[anchor=north west] {$s'_1$};
				\draw (1.35,2.75) node[anchor=north west] {$d'_1$};
				\draw (2.4,0.7) node[anchor=north west] {$s'_2$};
				\draw (1.35,0.75) node[anchor=north west] {$d'_2$};
				\draw (3.4,1.7) node[anchor=north west] {$s'_3$};
				\draw (2.35,1.75) node[anchor=north west] {$d'_3$};
				\draw (4.4,0.7) node[anchor=north west] {$s'_4$};
				\draw (3.35,0.75) node[anchor=north west] {$d'_4$};
				\draw (4.4,2.7) node[anchor=north west] {$s'_5$};
				
				\draw (5.4,1.7) node[anchor=north west] {$s'_1$};
				\draw (5.35,2.75) node[anchor=north west] {$d'_1$};
				\draw (6.4,0.7) node[anchor=north west] {$s'_2$};
				\draw (5.35,0.75) node[anchor=north west] {$d'_2$};
				\draw (7.4,1.7) node[anchor=north west] {$s'_3$};
				\draw (8.35,1.75) node[anchor=north west] {$d'_3$};
				\draw (6.4,2.7) node[anchor=north west] {$s'_4$};
				\draw (2.9,-0.1) node[anchor=north west] { $(i)$};
				\draw (6.9,-0.1) node[anchor=north west] { $(j)$};
				\begin{scriptsize}
					\draw [fill=red] (2.,2.) circle (2.5pt);
					\draw [fill=white] (3.,2.) circle (2.5pt);
					\draw [fill=white] (4.,2.) circle (2.5pt);
					\draw [fill=blue] (5.,2.) circle (2.5pt);
					\draw [fill=blue] (2.,1.) circle (2.5pt);
					\draw [fill=red] (3.,1.) circle (2.5pt);
					\draw [fill=blue] (4.,1.) circle (2.5pt);
					\draw [fill=white] (5.,1.) circle (2.5pt);
					\draw [fill=red] (2.,0.) circle (2.5pt);
					\draw [fill=blue] (3.,0.) circle (2.5pt);
					\draw [fill=red] (4.,0.) circle (2.5pt);
					\draw [fill=blue] (5.,0.) circle (2.5pt);
					\draw [fill=red] (6.,2.) circle (2.5pt);
					\draw [fill=blue] (7.,2.) circle (2.5pt);
					\draw [fill=white] (8.,2.) circle (2.5pt);
					\draw [fill=white] (9.,2.) circle (2.5pt);
					\draw [fill=blue] (6.,1.) circle (2.5pt);
					\draw [fill=white] (7.,1.) circle (2.5pt);
					\draw [fill=blue] (8.,1.) circle (2.5pt);
					\draw [fill=red] (9.,1.) circle (2.5pt);
					\draw [fill=red] (6.,0.) circle (2.5pt);
					\draw [fill=blue] (7.,0.) circle (2.5pt);
					\draw [fill=white] (8.,0.) circle (2.5pt);
					\draw [fill=white] (9.,0.) circle (2.5pt);
				\end{scriptsize}
				
				\draw [line width=1.2pt] (10.,2.)-- (11.,2.);
				\draw [line width=1.2pt] (11.,2.)-- (12.,2.);
				\draw [line width=1.2pt] (12.,2.)-- (13.,2.);
				\draw [line width=1.2pt] (10.,1.)-- (11.,1.);
				\draw [line width=1.2pt] (11.,1.)-- (12.,1.);
				\draw [line width=1.2pt] (12.,1.)-- (13.,1.);
				\draw [line width=1.2pt] (10.,0.)-- (11.,0.);
				\draw [line width=1.2pt] (11.,0.)-- (12.,0.);
				\draw [line width=1.2pt] (12.,0.)-- (13.,0.);
				\draw [line width=1.2pt] (10.,2.)-- (10.,1.);
				\draw [line width=1.2pt] (10.,1.)-- (10.,0.);
				\draw [line width=1.2pt] (11.,2.)-- (11.,1.);
				\draw [line width=1.2pt] (11.,1.)-- (11.,0.);
				\draw [line width=1.2pt] (12.,2.)-- (12.,1.);
				\draw [line width=1.2pt] (12.,1.)-- (12.,0.);
				\draw [line width=1.2pt] (13.,2.)-- (13.,1.);
				\draw [line width=1.2pt] (13.,1.)-- (13.,0.);
				\draw [line width=1.2pt] (14.,2.)-- (15.,2.);
				\draw [line width=1.2pt] (15.,2.)-- (16.,2.);
				\draw [line width=1.2pt] (16.,2.)-- (17.,2.);
				\draw [line width=1.2pt] (14.,1.)-- (15.,1.);
				\draw [line width=1.2pt] (15.,1.)-- (16.,1.);
				\draw [line width=1.2pt] (16.,1.)-- (17.,1.);
				\draw [line width=1.2pt] (14.,0.)-- (15.,0.);
				\draw [line width=1.2pt] (15.,0.)-- (16.,0.);
				\draw [line width=1.2pt] (16.,0.)-- (17.,0.);
				\draw [line width=1.2pt] (14.,2.)-- (14.,1.);
				\draw [line width=1.2pt] (14.,1.)-- (14.,0.);
				\draw [line width=1.2pt] (15.,2.)-- (15.,1.);
				\draw [line width=1.2pt] (15.,1.)-- (15.,0.);
				\draw [line width=1.2pt] (16.,2.)-- (16.,1.);
				\draw [line width=1.2pt] (16.,1.)-- (16.,0.);
				\draw [line width=1.2pt] (17.,2.)-- (17.,1.);
				\draw [line width=1.2pt] (17.,1.)-- (17.,0.);
				\draw (9.35,1.7) node[anchor=north west] {$s'_1$};
				\draw (9.35,2.75) node[anchor=north west] {$d'_1$};
				\draw (10.4,0.7) node[anchor=north west] {$s'_2$};
				\draw (9.35,0.75) node[anchor=north west] {$d'_2$};
				\draw (11.4,1.7) node[anchor=north west] {$s'_3$};
				\draw (11.4,0.75) node[anchor=north west] {$d'_3$};
				\draw (12.4,2.7) node[anchor=north west] {$s'_4$};
				
				\draw (13.4,1.7) node[anchor=north west] {$s'_1$};
				\draw (13.35,2.75) node[anchor=north west] {$d'_1$};
				\draw (14.4,0.7) node[anchor=north west] {$s'_2$};
				\draw (13.35,0.75) node[anchor=north west] {$d'_2$};
				\draw (15.4,1.7) node[anchor=north west] {$s'_3$};
				\draw (16.35,0.75) node[anchor=north west] {$d'_3$};
				\draw (15.4,2.7) node[anchor=north west] {$s'_4$};
				\draw (10.9,-0.1) node[anchor=north west] {$(k)$};
				\draw (14.9,-0.1) node[anchor=north west] {$(l)$};
				\begin{scriptsize}
					\draw [fill=red] (10.,2.) circle (2.5pt);
					\draw [fill=white] (11.,2.) circle (2.5pt);
					\draw [fill=white] (12.,2.) circle (2.5pt);
					\draw [fill=blue] (13.,2.) circle (2.5pt);
					\draw [fill=blue] (10.,1.) circle (2.5pt);
					\draw [fill=white] (11.,1.) circle (2.5pt);
					\draw [fill=blue] (12.,1.) circle (2.5pt);
					\draw [fill=white] (13.,1.) circle (2.5pt);
					\draw [fill=red] (10.,0.) circle (2.5pt);
					\draw [fill=blue] (11.,0.) circle (2.5pt);
					\draw [fill=red] (12.,0.) circle (2.5pt);
					\draw [fill=white] (13.,0.) circle (2.5pt);
					\draw [fill=red] (14.,2.) circle (2.5pt);
					\draw [fill=white] (15.,2.) circle (2.5pt);
					\draw [fill=blue] (16.,2.) circle (2.5pt);
					\draw [fill=white] (17.,2.) circle (2.5pt);
					\draw [fill=blue] (14.,1.) circle (2.5pt);
					\draw [fill=white] (15.,1.) circle (2.5pt);
					\draw [fill=blue] (16.,1.) circle (2.5pt);
					\draw [fill=white] (17.,1.) circle (2.5pt);
					\draw [fill=red] (14.,0.) circle (2.5pt);
					\draw [fill=blue] (15.,0.) circle (2.5pt);
					\draw [fill=white] (16.,0.) circle (2.5pt);
					\draw [fill=red] (17.,0.) circle (2.5pt);	
				\end{scriptsize}
			\end{tikzpicture}	
		\end{center}
		\begin{center}	
			\begin{tikzpicture}[line cap=round,line join=round,>=triangle 45,x=0.9cm,y=0.9cm]
				\draw [line width=1.2pt] (2.,2.)-- (3.,2.);
				\draw [line width=1.2pt] (3.,2.)-- (4.,2.);
				\draw [line width=1.2pt] (4.,2.)-- (5.,2.);
				\draw [line width=1.2pt] (2.,1.)-- (3.,1.);
				\draw [line width=1.2pt] (3.,1.)-- (4.,1.);
				\draw [line width=1.2pt] (4.,1.)-- (5.,1.);
				\draw [line width=1.2pt] (2.,0.)-- (3.,0.);
				\draw [line width=1.2pt] (3.,0.)-- (4.,0.);
				\draw [line width=1.2pt] (4.,0.)-- (5.,0.);
				\draw [line width=1.2pt] (2.,2.)-- (2.,1.);
				\draw [line width=1.2pt] (2.,1.)-- (2.,0.);
				\draw [line width=1.2pt] (3.,2.)-- (3.,1.);
				\draw [line width=1.2pt] (3.,1.)-- (3.,0.);
				\draw [line width=1.2pt] (4.,2.)-- (4.,1.);
				\draw [line width=1.2pt] (4.,1.)-- (4.,0.);
				\draw [line width=1.2pt] (5.,2.)-- (5.,1.);
				\draw [line width=1.2pt] (5.,1.)-- (5.,0.);
				\draw [line width=1.2pt] (6.,2.)-- (7.,2.);
				\draw [line width=1.2pt] (7.,2.)-- (8.,2.);
				\draw [line width=1.2pt] (8.,2.)-- (9.,2.);
				\draw [line width=1.2pt] (6.,1.)-- (7.,1.);
				\draw [line width=1.2pt] (7.,1.)-- (8.,1.);
				\draw [line width=1.2pt] (8.,1.)-- (9.,1.);
				\draw [line width=1.2pt] (6.,0.)-- (7.,0.);
				\draw [line width=1.2pt] (7.,0.)-- (8.,0.);
				\draw [line width=1.2pt] (8.,0.)-- (9.,0.);
				\draw [line width=1.2pt] (6.,2.)-- (6.,1.);
				\draw [line width=1.2pt] (6.,1.)-- (6.,0.);
				\draw [line width=1.2pt] (7.,2.)-- (7.,1.);
				\draw [line width=1.2pt] (7.,1.)-- (7.,0.);
				\draw [line width=1.2pt] (8.,2.)-- (8.,1.);
				\draw [line width=1.2pt] (8.,1.)-- (8.,0.);
				\draw [line width=1.2pt] (9.,2.)-- (9.,1.);
				\draw [line width=1.2pt] (9.,1.)-- (9.,0.);
				\draw (2.35,2.75) node[anchor=north west] {$d'_1$};
				\draw (4.4,1.7) node[anchor=north west] {$s'_2$};
				\draw (1.4,1.7) node[anchor=north west] {$s'_1$};
				\draw (4.4,2.7) node[anchor=north west] {$s'_3$};
				\draw (3.35,2.75) node[anchor=north west] {$d'_3$};
				\draw (4.4,0.7) node[anchor=north west] {$s'_4$};
				
				\draw (6.35,2.75) node[anchor=north west] {$d'_1$};
				\draw (5.4,1.7) node[anchor=north west] {$s'_1$};
				\draw (7.35,2.75) node[anchor=north west] {$d'_2$};
				\draw (8.4,1.7) node[anchor=north west] {$s'_2$};
				\draw (6.35,0.75) node[anchor=north west] {$d'_3$};
				\draw (5.4,0.7) node[anchor=north west] {$s'_3$};
				\draw (8.4,0.7) node[anchor=north west] {$s'_4$};
				\draw (2.9,-0.1) node[anchor=north west] { $(m)$};
				\draw (6.9,-0.1) node[anchor=north west] { $(n)$};
				\begin{scriptsize}
					\draw [color=red,fill=white] (2.,2.) circle (2.5pt);
					\draw [fill=red] (3.,2.) circle (2.5pt);
					\draw [fill=red] (4.,2.) circle (2.5pt);
					\draw [fill=blue] (5.,2.) circle (2.5pt);
					\draw [fill=blue] (2.,1.) circle (2.5pt);
					\draw [color=red, fill=white] (3.,1.) circle (2.5pt);
					\draw [fill=white] (4.,1.) circle (2.5pt);
					\draw [fill=blue] (5.,1.) circle (2.5pt);
					\draw [color=red,fill=white] (2.,0.) circle (2.5pt);
					\draw [color=red, fill=white] (3.,0.) circle (2.5pt);
					\draw [fill=white] (4.,0.) circle (2.5pt);
					\draw [fill=blue] (5.,0.) circle (2.5pt);
					\draw [fill=white] (6.,2.) circle (2.5pt);
					\draw [fill=red] (7.,2.) circle (2.5pt);
					\draw [fill=red] (8.,2.) circle (2.5pt);
					\draw [fill=white] (9.,2.) circle (2.5pt);
					\draw [fill=blue] (6.,1.) circle (2.5pt);
					\draw [fill=white] (7.,1.) circle (2.5pt);
					\draw [fill=white] (8.,1.) circle (2.5pt);
					\draw [fill=blue] (9.,1.) circle (2.5pt);
					\draw [fill=blue] (6.,0.) circle (2.5pt);
					\draw [fill=red] (7.,0.) circle (2.5pt);
					\draw [fill=white] (8.,0.) circle (2.5pt);
					\draw [fill=blue] (9.,0.) circle (2.5pt);
				\end{scriptsize}
				
				\draw [line width=1.2pt] (10.,2.)-- (11.,2.);
				\draw [line width=1.2pt] (11.,2.)-- (12.,2.);
				\draw [line width=1.2pt] (12.,2.)-- (13.,2.);
				\draw [line width=1.2pt] (10.,1.)-- (11.,1.);
				\draw [line width=1.2pt] (11.,1.)-- (12.,1.);
				\draw [line width=1.2pt] (12.,1.)-- (13.,1.);
				\draw [line width=1.2pt] (10.,0.)-- (11.,0.);
				\draw [line width=1.2pt] (11.,0.)-- (12.,0.);
				\draw [line width=1.2pt] (12.,0.)-- (13.,0.);
				\draw [line width=1.2pt] (10.,2.)-- (10.,1.);
				\draw [line width=1.2pt] (10.,1.)-- (10.,0.);
				\draw [line width=1.2pt] (11.,2.)-- (11.,1.);
				\draw [line width=1.2pt] (11.,1.)-- (11.,0.);
				\draw [line width=1.2pt] (12.,2.)-- (12.,1.);
				\draw [line width=1.2pt] (12.,1.)-- (12.,0.);
				\draw [line width=1.2pt] (13.,2.)-- (13.,1.);
				\draw [line width=1.2pt] (13.,1.)-- (13.,0.);
				\draw [line width=1.2pt] (14.,2.)-- (15.,2.);
				\draw [line width=1.2pt] (15.,2.)-- (16.,2.);
				\draw [line width=1.2pt] (16.,2.)-- (17.,2.);
				\draw [line width=1.2pt] (14.,1.)-- (15.,1.);
				\draw [line width=1.2pt] (15.,1.)-- (16.,1.);
				\draw [line width=1.2pt] (16.,1.)-- (17.,1.);
				\draw [line width=1.2pt] (14.,0.)-- (15.,0.);
				\draw [line width=1.2pt] (15.,0.)-- (16.,0.);
				\draw [line width=1.2pt] (16.,0.)-- (17.,0.);
				\draw [line width=1.2pt] (14.,2.)-- (14.,1.);
				\draw [line width=1.2pt] (14.,1.)-- (14.,0.);
				\draw [line width=1.2pt] (15.,2.)-- (15.,1.);
				\draw [line width=1.2pt] (15.,1.)-- (15.,0.);
				\draw [line width=1.2pt] (16.,2.)-- (16.,1.);
				\draw [line width=1.2pt] (16.,1.)-- (16.,0.);
				\draw [line width=1.2pt] (17.,2.)-- (17.,1.);
				\draw [line width=1.2pt] (17.,1.)-- (17.,0.);
				\draw (10.35,2.75) node[anchor=north west] {$d'_1$};
				\draw (12.4,1.7) node[anchor=north west] {$s'_2$};
				\draw (12.35,2.75) node[anchor=north west] {$d'_2$};
				\draw (9.4,1.7) node[anchor=north west] {$s'_1$};
				\draw (9.4,0.7) node[anchor=north west] {$s'_3$};
				\draw (10.35,0.75) node[anchor=north west] {$d'_3$};
				\draw (10.4,1.7) node[anchor=north west] {$s'_4$};
				\draw (9.35,2.75) node[anchor=north west] {$d'_4$};
				\draw (11.4,0.7) node[anchor=north west] {$s'_5$};
				
				\draw (14.35,2.75) node[anchor=north west] {$d'_1$};
				\draw (16.4,1.7) node[anchor=north west] {$s'_2$};
				\draw (15.35,1.75) node[anchor=north west] {$d'_2$};
				\draw (13.4,1.7) node[anchor=north west] {$s'_1$};
				\draw (16.4,2.7) node[anchor=north west] {$s'_3$};
				\draw (15.35,2.75) node[anchor=north west] {$d'_3$};
				\draw (15.4,0.7) node[anchor=north west] {$s'_4$};
				\draw (16.35,0.75) node[anchor=north west] {$d'_4$};
				\draw (13.4,0.7) node[anchor=north west] {$s'_5$};
				\draw (10.9,-0.1) node[anchor=north west] {$(o)$};
				\draw (14.9,-0.1) node[anchor=north west] {$(p)$};
				\begin{scriptsize}
					\draw [fill=red] (10.,2.) circle (2.5pt);
					\draw [fill=red] (11.,2.) circle (2.5pt);
					\draw [fill=white] (12.,2.) circle (2.5pt);
					\draw [fill=red] (13.,2.) circle (2.5pt);
					\draw [fill=blue] (10.,1.) circle (2.5pt);
					\draw [fill=blue] (11.,1.) circle (2.5pt);
					\draw [fill=white] (12.,1.) circle (2.5pt);
					\draw [fill=blue] (13.,1.) circle (2.5pt);
					\draw [fill=blue] (10.,0.) circle (2.5pt);
					\draw [fill=red] (11.,0.) circle (2.5pt);
					\draw [fill=blue] (12.,0.) circle (2.5pt);
					\draw [fill=white] (13.,0.) circle (2.5pt);
					\draw [fill=white] (14.,2.) circle (2.5pt);
					\draw [fill=red] (15.,2.) circle (2.5pt);
					\draw [fill=red] (16.,2.) circle (2.5pt);
					\draw [fill=blue] (17.,2.) circle (2.5pt);
					\draw [fill=blue] (14.,1.) circle (2.5pt);
					\draw [fill=white] (15.,1.) circle (2.5pt);
					\draw [fill=red] (16.,1.) circle (2.5pt);
					\draw [fill=blue] (17.,1.) circle (2.5pt);
					\draw [fill=blue] (14.,0.) circle (2.5pt);
					\draw [fill=white] (15.,0.) circle (2.5pt);
					\draw [fill=blue] (16.,0.) circle (2.5pt);
					\draw [fill=red] (17.,0.) circle (2.5pt);	
				\end{scriptsize}
			\end{tikzpicture}
		\end{center}
		
		\begin{center}
			\begin{tikzpicture}[line cap=round,line join=round,>=triangle 45,x=0.9cm,y=0.9cm]
				\draw [line width=1.2pt] (2.,2.)-- (3.,2.);
				\draw [line width=1.2pt] (3.,2.)-- (4.,2.);
				\draw [line width=1.2pt] (4.,2.)-- (5.,2.);
				\draw [line width=1.2pt] (2.,1.)-- (3.,1.);
				\draw [line width=1.2pt] (3.,1.)-- (4.,1.);
				\draw [line width=1.2pt] (4.,1.)-- (5.,1.);
				\draw [line width=1.2pt] (2.,0.)-- (3.,0.);
				\draw [line width=1.2pt] (3.,0.)-- (4.,0.);
				\draw [line width=1.2pt] (4.,0.)-- (5.,0.);
				\draw [line width=1.2pt] (2.,2.)-- (2.,1.);
				\draw [line width=1.2pt] (2.,1.)-- (2.,0.);
				\draw [line width=1.2pt] (3.,2.)-- (3.,1.);
				\draw [line width=1.2pt] (3.,1.)-- (3.,0.);
				\draw [line width=1.2pt] (4.,2.)-- (4.,1.);
				\draw [line width=1.2pt] (4.,1.)-- (4.,0.);
				\draw [line width=1.2pt] (5.,2.)-- (5.,1.);
				\draw [line width=1.2pt] (5.,1.)-- (5.,0.);
				\draw [line width=1.2pt] (6.,2.)-- (7.,2.);
				\draw [line width=1.2pt] (7.,2.)-- (8.,2.);
				\draw [line width=1.2pt] (8.,2.)-- (9.,2.);
				\draw [line width=1.2pt] (6.,1.)-- (7.,1.);
				\draw [line width=1.2pt] (7.,1.)-- (8.,1.);
				\draw [line width=1.2pt] (8.,1.)-- (9.,1.);
				\draw [line width=1.2pt] (6.,0.)-- (7.,0.);
				\draw [line width=1.2pt] (7.,0.)-- (8.,0.);
				\draw [line width=1.2pt] (8.,0.)-- (9.,0.);
				\draw [line width=1.2pt] (6.,2.)-- (6.,1.);
				\draw [line width=1.2pt] (6.,1.)-- (6.,0.);
				\draw [line width=1.2pt] (7.,2.)-- (7.,1.);
				\draw [line width=1.2pt] (7.,1.)-- (7.,0.);
				\draw [line width=1.2pt] (8.,2.)-- (8.,1.);
				\draw [line width=1.2pt] (8.,1.)-- (8.,0.);
				\draw [line width=1.2pt] (9.,2.)-- (9.,1.);
				\draw [line width=1.2pt] (9.,1.)-- (9.,0.);
				
				\draw (1.4,1.7) node[anchor=north west] {$s'_1$};
				\draw (2.35,2.75) node[anchor=north west] {$d'_1$};
				\draw (4.4,1.7) node[anchor=north west] {$s'_2$};
				\draw (3.35,0.75) node[anchor=north west] {$d'_2$};
				\draw (4.4,2.7) node[anchor=north west] {$s'_3$};
				\draw (3.35,2.75) node[anchor=north west] {$d'_3$};
				\draw (1.4,0.7) node[anchor=north west] {$s'_4$};	
				\draw (2.35,0.75) node[anchor=north west] {$d'_4$};

				\draw (5.4,0.7) node[anchor=north west] {$s'_3$};
				\draw (6.35,0.75) node[anchor=north west] {$d'_3$};
				\draw (5.4,1.7) node[anchor=north west] {$s'_1$};
				\draw (6.35,2.75) node[anchor=north west] {$d'_1$};
				\draw (8.4,1.7) node[anchor=north west] {$s'_2$};
				\draw (8.35,0.75) node[anchor=north west] {$d'_2$};
				\draw (7.4,2.7) node[anchor=north west] {$s'_4$};
				\draw (8.35,2.75) node[anchor=north west] {$d'_4$};
				\draw (6.4,1.7) node[anchor=north west] {$s'_5$};
				\draw (2.9,-0.1) node[anchor=north west] { $(q)$};
				\draw (6.9,-0.1) node[anchor=north west] { $(r)$};
				\begin{scriptsize}
					\draw [fill=white] (2.,2.) circle (2.5pt);
					\draw [fill=red] (3.,2.) circle (2.5pt);
					\draw [fill=red] (4.,2.) circle (2.5pt);
					\draw [fill=blue] (5.,2.) circle (2.5pt);
					\draw [fill=blue] (2.,1.) circle (2.5pt);
					\draw [fill=white] (3.,1.) circle (2.5pt);
					\draw [fill=white] (4.,1.) circle (2.5pt);
					\draw [fill=blue] (5.,1.) circle (2.5pt);
					\draw [fill=blue] (2.,0.) circle (2.5pt);
					\draw [fill=red] (3.,0.) circle (2.5pt);
					\draw [fill=red] (4.,0.) circle (2.5pt);
					\draw [fill=white] (5.,0.) circle (2.5pt);
					\draw [fill=white] (6.,2.) circle (2.5pt);
					\draw [fill=red] (7.,2.) circle (2.5pt);
					\draw [fill=blue] (8.,2.) circle (2.5pt);
					\draw [fill=red] (9.,2.) circle (2.5pt);
					\draw [fill=blue] (6.,1.) circle (2.5pt);
					\draw [fill=blue] (7.,1.) circle (2.5pt);
					\draw [fill=white] (8.,1.) circle (2.5pt);
					\draw [fill=blue] (9.,1.) circle (2.5pt);
					\draw [fill=blue] (6.,0.) circle (2.5pt);
					\draw [fill=red] (7.,0.) circle (2.5pt);
					\draw [fill=white] (8.,0.) circle (2.5pt);
					\draw [fill=red] (9.,0.) circle (2.5pt);
				\end{scriptsize}
				
			\end{tikzpicture}
			\caption{An illustration for Staller's strategy in the proof of $\gmb'(P_3 \cp P_4) \ge 6$. Blue and red vertices denote the moves of Staller and Dominator, respectively, and vertices which are outlined in red denote possible moves of Dominator.}
			\label{fig:Stall Sgame}
		\end{center}
	\end{figure}

	\paragraph{Case 1.} Dominator plays $d'_1 \in [3]\times \{3,4\}$. Then Staller responds by playing $s'_2=(3,1)$. Dominator needs to reply $d'_2=(3,2)$, otherwise Staller will win in the next move. After that Staller plays $s'_3=(1,1)$, Dominator needs to play $d'_3 \in \{(1,2), (2,2)\}$.
	Then Staller will win in her next move by playing the unplayed vertex from the set. See Figure \ref{fig:Stall Sgame}\,(a).
	
	\paragraph{Case 2.} Dominator plays $d'_1 =(2,2)$. Then Staller replies by playing $s'_2=(1,2)$ and it forces Dominator to play $d'_2 =(1,1)$. Then Staller plays $s'_3=(3,2)$ and Dominator must play $d'_3 =(3,1)$. After that Staller plays $s'_4=(3,4)$. 
	
	\begin{itemize}
		\item[(2.1)] If $d'_4 =(3,3)$, then Staller responds by selecting $s'_5=(1,4)$. 
		\item[(2.2)] If $d'_4 =(2,3)$, then Staller responds by selecting $s'_5=(2,4)$.
		\item[(2.3)] If $d'_4 =(2,4)$, then Staller responds by selecting $s'_5=(2,3)$.
		\item[(2.4)] If $d'_4 =(1,3)$, then Staller responds by selecting $s'_5=(2,4)$.
		\item[(2.5)] If $d'_4 =(1,4)$, then Staller responds by selecting $s'_5=(3,3)$.	
	\end{itemize}
	By above strategy, Dominator cannot dominate $Z$ with five moves if his first move is $d'_1=(2,2)$. See Figure \ref{fig:Stall Sgame}\,(b)-(f).
	
	\paragraph{Case 3.} Dominator plays $d'_1 =(3,1)$. Then Staller plays $s'_2=(1,2)$ and Dominator must select $d'_2=(1,1)$. After that Staller plays $s'_3=(2,3)$.
	\begin{itemize}
		\item[(3.1)] If $d'_3 =(3,2)$, then Staller selects $s'_4=(2,4)$. 
		If $d'_4=(1,3)$, then Staller plays $s'_5=(3,4)$. 
		If $d'_4=(1,4)$, then Staller plays $s'_5=(3,3)$.
		If $d'_4=(3,3)$, then Staller plays $s'_5=(1,4)$. 
		If $d'_4=(3,4)$, then Staller plays $s'_5=(1,3)$. In each case, Dominator cannot dominate $Z$ within five moves.
		\item[(3.2)] If $d'_3 =(3,3)$, then Staller responds by selecting $s'_4=(1,4)$.	
		If $d'_4 \neq (1,3)$, then Staller plays $s'_5=(1,3)$ and wins the game. So, Dominator needs to play $d'_4 = (1,3)$ and he also needs two more moves to dominate $Z$. Thus he needs to play at least six moves.
		\item[(3.3)] If $d'_3 =(2,2)$, then Staller responds by selecting $s'_4=(1,4)$ and Dominator has to play $d'_4=(1,3)$. Thus Staller plays $s'_5=(3,4)$ and Dominator needs two more moves to win the game.
		\item[(3.4)] If $d'_3 =(2,4)$, then Staller responds by selecting $s'_4=(3,2)$. Thus Dominator needs three more moves to dominate $Z$.
		\item[(3.5)] If $d'_3 =(1,3)$, then Staller responds by selecting $s'_4=(3,4)$. Then Staller plays $s'_5=(2,4)$ if it is possible, otherwise $s'_5=(3,2)$. Thus Dominator needs to play at least six moves.
		\item[(3.6)] If $d'_3 =(1,4)$, then Staller responds by selecting $s'_4=(3,3)$. Observe that only $(2,2), (2,3), (3,3), (3,4)$ are undominated and $(2,3), (3,3)$ are played by Staller. 
		If $d'_4=(2,2)$, then Staller plays $s'_5=(3,4)$.
		If $d'_4=(2,4)$, then Staller plays $s'_5=(3,2)$.
		If $d'_4=(3,2)$, then Staller plays $s'_5=(2,4)$. 
		If $d'_4=(3,4)$, then Staller plays $s'_5=(2,2)$. 
		In each case, Dominator cannot dominate $Z$ within five moves.
	\end{itemize}
	Therefore, Dominator needs at least six moves to win the game if he starts the game with $d'_1=(3,1)$. See Figure \ref{fig:Stall Sgame}\,(g)-(l).
	\paragraph{Case 4.} Dominator plays $d'_1 =(3,2)$. Then Staller replies by playing $s'_2=(2,4)$.
	\begin{itemize}
		\item If $d'_2 \in \{ (1,1), (1,2), (2,2), (3,1)\}$, then Staller replies by playing $s'_3=(3,4)$ and Dominator has to play $d'_3=(3,3)$, otherwise Staller will win in the next move. After that Staller plays $s'_4=(1,4)$ and she will win in her next turn.
		\item If $d'_2= (3,3)$, then Staller replies by playing $s'_3=(1,1)$ and Dominator has to play $d'_3=(1,2)$, otherwise Staller will win in the next move. Later Staller plays $s'_4=(1,4)$ and Dominator needs three more moves to win the game.
		\item If $d'_2= (3,4)$, then Staller replies by playing $s'_3=(1,1)$ and Dominator has to play $d'_3=(1,2)$, otherwise Staller will win in the next move. Later Staller plays $s'_4=(2,2)$ and it forces Dominator to play $d'_4=(3,1)$. Thus Dominator needs two more moves to win the game.
		\item If $d'_2= (2,3)$, then Staller replies by playing $s'_3=(3,4)$ and Dominator needs to play $d'_3=(3,3)$. After that Staller plays $s'_4=(1,4)$ an it forces Dominator to play $d'_4=(1,3)$. Thus Staller plays $s'_5=(1,1)$ and Dominator needs two more moves to win the game.
		\item If $d'_2= (1,3)$, then Staller replies by playing $s'_3=(3,4)$ and Dominator needs to play $d'_3=(3,3)$. After that Staller plays $s'_4=(1,1)$ and Dominator needs to play $d'_4=(1,2)$. Thus Dominator needs to play two more moves to win the game.
		\item If $d'_2= (1,4)$, then Staller replies by playing $s'_3=(1,1)$ and it forces Dominator to reply $d'_3=(1,2)$. Then Staller plays $s'_4=(3,3)$ and it forces Dominator to play $d'_4=(3,4)$. Later Staller plays $s'_5=(2,2)$.
		Thus Dominator need two more moves to win the game.
	\end{itemize}
	By the above strategy, Dominator needs at least six moves to win the game if his first move is $d'_1=(3,2)$. See Figure \ref{fig:Stall Sgame}\,(m)-(r).
	
	\medskip
	These four cases, together with the symmetrical ones, cover all possibilities for Dominator's first move. We conclude that $\gmb' (Z) \ge 6$.
		\qed
		
	\medskip
		In~\cite{jacobson-1984}, the domination number of $P_m \cp P_n$ was established for $1 \le m \le 4$ and all $n \ge 1$. In particular, it was shown that \begin{align} \label{gam_grid}
		\gamma(P_3 \cp P_n) = \left\lfloor\frac{3n+4}{4}\right\rfloor,\ n \ge 1.
	\end{align}
	
	\begin{theorem} \label{thm:boundP_3P_n}
		If $n \ge 2$, then $\lfloor\frac{3n+4}{4}\rfloor \le \gmb(P_3 \cp P_n) \le \frac{4n +\sigma(n)}{3}$ where 
		$$\sigma(n)= \begin{cases} 
				0; & 3|n,\\
				-1; & n \equiv 1 \pmod {3},\\ 
				1; &  n \equiv 2 \pmod {3}.\\
				\end{cases}$$ 
	\end{theorem}
	\proof By~(\ref{gam_grid}), $\gmb(P_3 \cp P_n) \ge \gamma(P_3 \cp P_n) = \lfloor\frac{3n+4}{4}\rfloor$. It remains to show that $\gmb(P_3 \cp P_n) \le \lceil\frac{4n}{3}\rceil$.
	It is easy to see that $\gmb(P_3 \cp P_2) = \gmb'(P_3 \cp P_2)=  3$.	For $n \ge 4$ we consider the following cases. 
	
	\paragraph{Case 1.} $n=3k$ where $k \ge 1$. Then $\sigma(n)=0$.
	Let $G$ be a disjoint union\index{disjoint union} of $k$ copies of $P_3 \cp P_3$. By Proposition~\ref{prop:p3p3} and Theorem~\ref{thm:gmb-union}\,, $\gmb(G) = 4k$. Since $G$ is a spanning subgraph of $P_3 \cp P_n$, Lemma~\ref{lem:o_edgedelete}\,(i) implies $\gmb(P_3 \cp P_n) \le 4k =\frac{4n}{3}$.
	
	\paragraph{Case 2.} $n=3k+1$ where $k \ge 1$. Then $\sigma(n)=-1$.
	Let $G$ be a union of $k-1$ copies of $P_3 \cp P_3$ and a copy of $P_3 \cp P_4$. By Theorem~\ref{thm:P2}, Proposition~\ref{prop:p3p3}, and Theorem~\ref{thm:gmb-union}\,, $\gmb(G) \le 4(k-1)+5$. Similar to the proof of Case 1, $\gmb(P_3 \cp P_n) \le 4k+1=\frac{4(n-1)}{3}+1 =\frac{4n-1}{3}$.
	\paragraph{Case 3.}$n=3k+2$  where $k \ge 1$. Then $\sigma(n)=1$.
	Let $G$ be a union of $k$ copies of $P_3 \cp P_3$ and a copy of $P_3 \cp P_2$. By Proposition~\ref{prop:p3p3} and Theorem~\ref{thm:gmb-union}\,, $\gmb(G) = 4k+3$. Similar to the proof of Case 1, $\gmb(P_3 \cp P_n) \le 4k+3=\frac{4(n-2)}{3}+3 =\frac{4n+1}{3}$.
	
	\medskip
	We conclude that  $\gmb(P_3 \cp P_n) \le \frac{4n +\sigma(n)}{3}$ for every $n \ge 2$.
	\qed
	\medskip

	\section{Products of complete bipartite graphs}\label{sec: product-bipart}
	In this section, we consider the MBD game on Cartesian products of complete bipartite graphs. We determine the outcome of the game, the exact MBD-numbers, and SMBD-numbers over this graph class.
	\begin{proposition}\label{prop:o_Kmn}
		If $n \ge m \ge 1$, then
		\begin{align*}
			o(K_{m,n}) &= \begin{cases}
				\cal D; & n=m=1,\text{or\ }n, m \ge 2, \\
				\cal N; & n > m =1.
			\end{cases}
		\end{align*}
	\end{proposition}
	\proof
	It is easy to see that Dominator wins the game in his first move on $K_{1,1}$ in both games. Thus $o(K_{1,1})=\cD.$
	
	Now consider star $K_{1,n}$, $n \ge 2$. In the D-game, Dominator plays the central vertex and wins the game in his first move. In the S-game, Staller also plays the central vertex in her first move. No matter where Dominator replies, Staller wins in her next move by playing an unplayed leaf. Thus $o(K_{1,n}) = \cal N$.
	
	Next, assume that $n \ge 2$ and $m \ge 2$. Let $X$, $Y$ be the bipartition sets of $K_{m,n}$. By Corollary \ref{cor:skip}, it suffices to show that Dominator can win the S-game. Without loss of generality, suppose that Staller starts the game by playing a vertex in $X$. Then Dominator plays an unplayed vertex in $X$ and dominates all vertices in $Y$. Then Dominator plays an unplayed vertex in $Y$ and wins the game.  We conclude that $o(K_{m,n}) = \cal D$, where $n, m \ge 2$.
	\qed
	
	Proposition~\ref{prop:o_Kmn} can also be deduced by combining the result from \cite{gledel-2020} asserting that every cograph which admits a paring dominating set has outcome $\cal D$, and the fact that complete bipartite graphs are cographs.
	
	By Theorem~\ref{thm:cart} and Proposition~\ref{prop:o_Kmn}\,, we obtain the outcomes of the MBD game on the Cartesian product of complete bipartite graphs as follows.
	\begin{corollary}
		Let $m, m', n, n'$ be positive integers. If $n \ge m \ge 2$ and $n' \ge m' \ge 2$, then $o(K_{m,n} \cp K_{m',n'}) = \cal D$.
	\end{corollary}
	
	Note that Theorem~\ref{thm:cart} does not help us to find the outcome of products of two stars. Moreover, a star $K_{1,n}$, where  $n \ge 3$ does not admit a nontrivial path cover by Theorem~\ref{thm: factor}\,. Hence, we need to consider the outcome of products of stars as follows.
	
	Let $Z$ be the graph $K_{1,m} \cp K_{1,n}$, and  $V(K_{1,m} \cp K_{1,n}) = \{(i,j):\ i\in [m], j\in [n]\} \cup \{(a,j):\ j\in [n]\} \cup \{(i,b):\ i\in [m]\}$, where $a, b$ are the central vertices of  $K_{1,m}$, and $K_{1,n}$, respectively.
	
	\begin{figure}[h!]
		\begin{center}	
			\begin{tikzpicture}[line cap=round,line join=round,>=triangle 45,x=0.9cm,y=0.9cm]
				\draw [line width=1.2pt] (-3.,2.)-- (-2.,2.);
				\draw [line width=1.2pt] (-3.,3.)-- (-2.,3.);
				\draw [line width=1.2pt] (-3.,4.)-- (-2.,4.);
				\draw [line width=1.2pt] (-3.,5.8)-- (-2.,5.8);
				\draw [line width=1.2pt] (-3.,2.)-- (-3.,3.);
				\draw [line width=1.2pt] (-2.,2.)-- (-2.,3.);
				\draw [line width=1.2pt] (-1.,2.)-- (-1.,3.);
				\draw [line width=1.2pt] (0.,2.)-- (0.,3.);
				\draw [line width=1.2pt] (2.,2.)-- (2.,3.);
				\draw [line width=1.2pt,color=black] (-3.,5.8)..controls(-3.,6.1)and(-1.,6.1).. (-1.,5.8);
				\draw [line width=1.2pt,color=black] (-3.,5.8)..controls(-3.,6.1)and(0,6.1).. (0.,5.8);
				\draw [line width=1.2pt,color=black] (-3.,4.)..controls(-3.,4.3)and(-1.,4.3).. (-1.,4.);
				\draw [line width=1.2pt,color=black] (-3.,4.)..controls(-3,4.3)and(0,4.3).. (0.,4.);
				\draw [line width=1.2pt,color=black] (-3.,3.)..controls(-3.,3.3)and(-1.,3.3)..  (-1.,3.);
				\draw [line width=1.2pt,color=black] (-3.,3.)..controls(-3,3.3)and(0,3.3).. (0.,3.);
				\draw [line width=1.2pt,color=black] (-3.,2.)..controls(-3.,2.3)and(-1.,2.3)..  (-1.,2.);
				\draw [line width=1.2pt,color=black] (-3.,2.)..controls(-3,2.3)and(0,2.3).. (0.,2.);
				\draw [line width=1.2pt,color=black] (2.,5.8)..controls(2,6.2)and (-3,6.2)..(-3.,5.8);
				\draw [line width=1.2pt,color=black]  (2.,4.)..controls(2,4.4)and (-3,4.4)..(-3.,4.);
				
				\draw [line width=1.2pt,color=black] (2.,3.)..controls(2,3.4)and (-3,3.4)..(-3.,3.);
				\draw [line width=1.2pt,color=black]  (2.,2.)..controls(2,2.4)and (-3,2.4)..(-3.,2.);

				\draw [line width=1.2pt,color=black] (-3.,2.)..controls(-3.3,2.)and (-3.3,4.)..(-3.,4.);
				\draw [line width=1.2pt,color=black] (-3.,2.)..controls(-3.3,2.)and (-3.3,5.8)..(-3.,5.8);
				\draw [line width=1.2pt,color=black] (-2.,2.)..controls(-2.3,2.)and (-2.3,4.)..(-2.,4.);
				\draw [line width=1.2pt,color=black] (-2.,2.)..controls(-2.3,2.)and (-2.3,5.8)..(-2.,5.8);
				\draw [line width=1.2pt,color=black] (-1.,2.)..controls(-1.3,2.)and (-1.3,4.)..(-1.,4.);
				\draw [line width=1.2pt,color=black] (-1.,2.)..controls(-1.3,2.)and (-1.3,5.8)..(-1.,5.8);
				\draw [line width=1.2pt,color=black] (0.,2.)..controls(-0.3,2.)and (-0.3,4.)..(0.,4.);
				\draw [line width=1.2pt,color=black] (0.,2.)..controls(-0.3,2.)and (-0.3,5.8)..(0.,5.8);
				\draw [line width=1.2pt,color=black] (2.,2.)..controls(1.7,2.)and (1.7,4.)..(2.,4.);
				\draw [line width=1.2pt,color=black] (2.,2.)..controls(1.7,2.)and (1.7,5.8)..(2.,5.8);

				\begin{scriptsize}
					\draw [fill=white] (-3.,5.8) circle (2.5pt);
					\draw (-3.25,1.85) node[anchor=north west] {\footnotesize$b$};
					\draw (-3.75,2.25) node[anchor=north west] {\footnotesize$a$};
					\draw [fill=white] (-3.,4.) circle (2.5pt);
					\draw (-3.75,3.25) node[anchor=north west] {\footnotesize$1$};
					\draw [fill=white] (-3.,3.) circle (2.5pt);
					\draw (-3.75,4.25) node[anchor=north west] {\footnotesize$2$};
					\draw [fill=white] (-3.,2.) circle (2.5pt);
					\draw (-3.75,6) node[anchor=north west] {\footnotesize$m$};
					\draw (-3.75,5.5) node[anchor=north west] {\footnotesize$\vdots$};
					\draw (-2.75,5.5) node[anchor=north west] {\footnotesize$\vdots$};
					\draw (-1.75,5.5) node[anchor=north west] {\footnotesize$\vdots$};
					\draw (-0.75,5.5) node[anchor=north west] {\footnotesize$\vdots$};
					\draw (0.55,5.5) node[anchor=north west] {\footnotesize$\iddots$};
					\draw [fill=white] (-2.,2.) circle (2.5pt);
					\draw (-2.25,1.85) node[anchor=north west] {\footnotesize$1$};
					\draw [fill=white] (-2.,3.) circle (2.5pt);
					\draw [fill=white] (-2.,4.) circle (2.5pt);
					\draw [fill=white] (-2.,5.8) circle (2.5pt);
					\draw [fill=white] (-1.,5.8) circle (2.5pt);
					\draw [fill=white] (-1.,4.) circle (2.5pt);
					\draw [fill=white] (-1.,3.) circle (2.5pt);
					\draw [fill=white] (-1.,2.) circle (2.5pt);
					\draw (-1.25,1.85) node[anchor=north west] {\footnotesize$2$};
					\draw [fill=white] (0.,2.) circle (2.5pt);
					\draw (-0.25,1.85) node[anchor=north west] {\footnotesize$3$};
					\draw [fill=white] (0.,3.) circle (2.5pt);
					\draw [fill=white] (0.,4.) circle (2.5pt);
					\draw [fill=white] (0.,5.8) circle (2.5pt);
					\draw [fill=white] (2.,2.) circle (2.5pt);
					\draw [fill=white] (2.,3.) circle (2.5pt);
					\draw [fill=white] (2.,4.) circle (2.5pt);
					\draw [fill=white] (2.,5.8) circle (2.5pt);
					\draw (0.65,1.85) node[anchor=north west] {\footnotesize$\ldots$};
					\draw (0.65,2.85) node[anchor=north west] {\footnotesize$\ldots$};
					\draw (0.65,3.85) node[anchor=north west] {\footnotesize$\ldots$};
					\draw (1.85,1.85) node[anchor=north west] {\footnotesize$n$};
				\end{scriptsize}
			\end{tikzpicture}
			\caption{graph $Z=K_{1,m} \cp K_{1,n}$.}
			\label{fig:Z}
		\end{center}
	\end{figure}
		
	\begin{theorem}\label{prop:out_star}
		If $n \ge m \ge 2$, then
		\begin{align}
			o(K_{1,m} \cp K_{1,n}) &= \begin{cases}
				\cal D; & m=n=2,\ \\
				\cal N; & m=2,\ n \ge 3,\\
				\cal S; & n \ge m \ge 3.
			\end{cases}
		\end{align}
	\end{theorem}
	\proof
	Assume that $n \ge m \ge 2$ and  $Z = K_{1,m} \cp K_{1,n}$ We consider the following cases.
	\paragraph{Case 1.} If $m=n=2$, then $Z=K_{1,2} \cp K_{1,2} = P_3 \cp P_3$. By Theorem~\ref{thm:pmpn}\,, $o(Z) = \cal D$.
	
	\paragraph{Case 2.} If $m=2$ and $n \ge 3$, we will show that the first player has a winning strategy on $Z = K_{1,2} \cp K_{1,n}$, $n \ge 3$ in the D-game and the S-game.
	
	In the D-game, Dominator starts the game by playing $(1,b)$. Then there is a matching $M$ in $G-\{(1,b)\}$ such that $V(Z)\setminus V(M) \subseteq N_Z[(1,b)]$. By Lemma~\ref{lem:pairing}\,, Dominator wins the game.
	
	In the S-game, Staller plays $s'_1=(1,b)$ in her first move. After first move $d'_1$ of Dominator, we consider the following strategies.
	\begin{itemize}
		\item $d'_1=(2,b)$. Then Staller continues the game by playing $s'_2=(1,1)$ that forces Dominator to play $d'_2=(a,1)$. Otherwise, Staller can win in the next move. In each turn, Staller uses the same strategy to play $s'_{i+1}=(1,i)$  for each $i \in [n]$ until she plays all the vertices in the layer $^1(K_{1,n})$, and she will win the game by playing $(a,b)$ in her last move.
		
		\item $d'_1 \neq (2,b)$. Then Staller responds on vertex $s'_2=(2,b)$. After the second move of Dominator, there is an unplayed layer of $(K_{1,2})^j$ and Staller will play $s'_3=(a,j)$. Thus Staller can win the game in her next move.
	\end{itemize}
	Thus $o( K_{1,2} \cp K_{1,n}) = \cal N$, where $n \ge 3.$
	\paragraph{Case 3.} Assume that $m=3$ and  $n = 3$. To show that $o(Z) = \cal S$, by Corollary \ref{cor:skip}\,(ii), it suffices to show that Staller has a winning strategy in the D-game when $Z =K_{1,3} \cp K_{1,3}$. Consider the following possible cases.
	\medskip

\begin{figure}[h]
		\begin{center}	
			\begin{tikzpicture}[line cap=round,line join=round,>=triangle 45,x=0.9cm,y=0.9cm]
				\clip(-4.,-4.) rectangle (11.02,5.88);
				\draw [line width=1.2pt] (-3.,2.)-- (-2.,2.);
				\draw [line width=1.2pt] (-3.,3.)-- (-2.,3.);
				\draw [line width=1.2pt] (-3.,4.)-- (-2.,4.);
				\draw [line width=1.2pt] (-3.,5.)-- (-2.,5.);
				\draw [line width=1.2pt] (-3.,2.)-- (-3.,3.);
				\draw [line width=1.2pt] (-2.,2.)-- (-2.,3.);
				\draw [line width=1.2pt] (-1.,2.)-- (-1.,3.);
				\draw [line width=1.2pt] (0.,2.)-- (0.,3.);
				\draw [line width=1.2pt] (2.,2.)-- (2.,3.);
				\draw [line width=1.2pt] (3.,2.)-- (3.,3.);
				\draw [line width=1.2pt] (4.,2.)-- (4.,3.);
				\draw [line width=1.2pt] (5.,2.)-- (5.,3.);
				\draw [line width=1.2pt] (2.,2.)-- (3.,2.);
				\draw [line width=1.2pt] (2.,3.)-- (3.,3.);
				\draw [line width=1.2pt] (2.,4.)-- (3.,4.);
				\draw [line width=1.2pt] (2.,5.)-- (3.,5.);
				\draw [line width=1.2pt] (7.,2.)-- (7.,3.);
				\draw [line width=1.2pt] (8.,2.)-- (8.,3.);
				\draw [line width=1.2pt] (9.,2.)-- (9.,3.);
				\draw [line width=1.2pt] (10.,2.)-- (10.,3.);
				\draw [line width=1.2pt] (7.,2.)-- (8.,2.);
				\draw [line width=1.2pt] (7.,3.)-- (8.,3.);
				\draw [line width=1.2pt] (7.,4.)-- (8.,4.);
				\draw [line width=1.2pt] (7.,5.)-- (8.,5.);
				\draw [line width=1.2pt] (-3.,-3.)-- (-2.,-3.);
				\draw [line width=1.2pt] (-3.,-2.)-- (-2.,-2.);
				\draw [line width=1.2pt] (-3.,-1.)-- (-2.,-1.);
				\draw [line width=1.2pt] (-3.,0.)-- (-2.,0.);
				\draw [line width=1.2pt] (-3.,-2.)-- (-3.,-3.);
				\draw [line width=1.2pt] (-2.,-2.)-- (-2.,-3.);
				\draw [line width=1.2pt] (-1.,-2.)-- (-1.,-3.);
				\draw [line width=1.2pt] (0.,-2.)-- (0.,-3.);
				\draw [line width=1.2pt] (2.,-3.)-- (3.,-3.);
				\draw [line width=1.2pt] (2.,-2.)-- (3.,-2.);
				\draw [line width=1.2pt] (2.,-1.)-- (3.,-1.);
				\draw [line width=1.2pt] (2.,0.)-- (3.,0.);
				\draw [line width=1.2pt] (2.,-3.)-- (2.,-2.);
				\draw [line width=1.2pt] (3.,-3.)-- (3.,-2.);
				\draw [line width=1.2pt] (4.,-3.)-- (4.,-2.);
				\draw [line width=1.2pt] (5.,-3.)-- (5.,-2.);
				\draw [line width=1.2pt] (7.,-2.)-- (8.,-2.);
				\draw [line width=1.2pt] (7.,-3.)-- (8.,-3.);
				\draw [line width=1.2pt] (7.,-1.)-- (8.,-1.);
				\draw [line width=1.2pt] (7.,0.)-- (8.,0.);
				\draw [line width=1.2pt] (7.,-3.)-- (7.,-2.);
				\draw [line width=1.2pt] (8.,-3.)-- (8.,-2.);
				\draw [line width=1.2pt] (9.,-3.)-- (9.,-2.);
				\draw [line width=1.2pt] (10.,-3.)-- (10.,-2.);
				\draw [line width=1.2pt,color=black] (-3.,5.)..controls(-3.,5.3)and(-1.,5.3).. (-1.,5.);
				\draw [line width=1.2pt,color=black] (-3.,5.)..controls(-3.,5.3)and(0,5.3).. (0.,5.);
				\draw [line width=1.2pt,color=black] (-3.,4.)..controls(-3.,4.3)and(-1.,4.3).. (-1.,4.);
				\draw [line width=1.2pt,color=black] (-3.,4.)..controls(-3,4.3)and(0,4.3).. (0.,4.);
				\draw [line width=1.2pt,color=black] (-3.,3.)..controls(-3.,3.3)and(-1.,3.3)..  (-1.,3.);
				\draw [line width=1.2pt,color=black] (-3.,3.)..controls(-3,3.3)and(0,3.3).. (0.,3.);
				\draw [line width=1.2pt,color=black] (-3.,2.)..controls(-3.,2.3)and(-1.,2.3)..  (-1.,2.);
				\draw [line width=1.2pt,color=black] (-3.,2.)..controls(-3,2.3)and(0,2.3).. (0.,2.);
				\draw [line width=1.2pt,color=black] (2.,5.)..controls(2,5.3)and (4,5.3)..(4.,5.);
				\draw [line width=1.2pt,color=black] (2.,5.)..controls(2,5.3)and (5,5.3)..(5.,5.);
				\draw [line width=1.2pt,color=black] (2.,4.)..controls(2,4.3)and (4,4.3)..(4.,4.);
				\draw [line width=1.2pt,color=black] (2.,4.)..controls(2,4.3)and (5,4.3)..(5.,4.);
				\draw [line width=1.2pt,color=black] (2.,3.)..controls(2,3.3)and (4,3.3)..(4.,3.);
				\draw [line width=1.2pt,color=black] (2.,3.)..controls(2,3.3)and (5,3.3)..(5.,3.);
				\draw [line width=1.2pt,color=black] (2.,2.)..controls(2,2.3)and (4,2.3)..(4.,2.);
				\draw [line width=1.2pt,color=black] (2.,2.)..controls(2,2.3)and (5,2.3)..(5.,2.);
				\draw [line width=1.2pt,color=black] (7.,5.)..controls(7,5.3)and (9,5.3)..(9.,5.);
				\draw [line width=1.2pt,color=black] (7.,5.)..controls(7,5.3)and (10,5.3)..(10.,5.);
				\draw [line width=1.2pt,color=black] (7.,4.)..controls(7,4.3)and (9,4.3)..(9.,4.);
				\draw [line width=1.2pt,color=black] (7.,4.)..controls(7,4.3)and (10,4.3)..(10.,4.);
				\draw [line width=1.2pt,color=black] (7.,3.)..controls(7.,3.3)and (9,3.3)..(9.,3.);
				\draw [line width=1.2pt,color=black] (7.,3.)..controls(7.,3.3)and (10,3.3)..(10.,3.);
				\draw [line width=1.2pt,color=black] (7.,2.)..controls(7.,2.3)and (9,2.3)..(9.,2.);
				\draw [line width=1.2pt,color=black] (7.,2.)..controls(7.,2.3)and (10,2.3)..(10.,2.);
				\draw [line width=1.2pt,color=black] (-3.,0.)..controls(-3.,0.3)and (-1,0.3)..(-1.,0.);
				\draw [line width=1.2pt,color=black] (-3.,0.)..controls(-3.,0.3)and (0,0.3)..(0.,0.);
				\draw [line width=1.2pt,color=black] (-3.,-1.)..controls(-3.,-0.7)and (-1,-0.7)..(-1.,-1.);
				\draw [line width=1.2pt,color=black] (-3.,-1.)..controls(-3.,-0.7)and (0,-0.7)..(0.,-1.);
				\draw [line width=1.2pt,color=black] (-3.,-2.)..controls(-3.,-1.7)and (-1,-1.7)..(-1.,-2.);
				\draw [line width=1.2pt,color=black] (-3.,-2.)..controls(-3.,-1.7)and (0,-1.7)..(0.,-2.);
				\draw [line width=1.2pt,color=black] (-3.,-3.)..controls(-3.,-2.7)and (-1,-2.7)..(-1.,-3.);
				\draw [line width=1.2pt,color=black] (-3.,-3.)..controls(-3.,-2.7)and (0.,-2.7)..(0.,-3.);
				\draw [line width=1.2pt,color=black] (2.,0.)..controls(2.,0.3)and (4.,0.3)..(4.,0.);
				\draw [line width=1.2pt,color=black] (2.,0.)..controls(2,0.3)and (5,0.3)..(5.,0.);
				\draw [line width=1.2pt,color=black] (2.,-1.)..controls(2,-0.7)and (4,-0.7)..(4.,-1.);
				\draw [line width=1.2pt,color=black] (2.,-1.)..controls(2,-0.7)and (5,-0.7)..(5.,-1.);
				\draw [line width=1.2pt,color=black] (2.,-2.)..controls(2,-1.7)and (4,-1.7)..(4.,-2.);
				\draw [line width=1.2pt,color=black] (2.,-2.)..controls(2,-1.7)and (5,-1.7)..(5.,-2.);
				\draw [line width=1.2pt,color=black] (2.,-3.)..controls(2,-2.7)and (4,-2.7)..(4.,-3.);
				\draw [line width=1.2pt,color=black] (2.,-3.)..controls(2,-2.7)and (5,-2.7)..(5.,-3.);
				\draw [line width=1.2pt,color=black] (7.,0.)..controls(7,0.3)and (9,0.3)..(9.,0.);
				\draw [line width=1.2pt,color=black] (7.,0.)..controls(7,0.3)and (10,0.3)..(10.,0.);
				\draw [line width=1.2pt,color=black] (7.,-1.)..controls(7,-0.7)and (9,-0.7)..(9.,-1.);
				\draw [line width=1.2pt,color=black] (7.,-1.)..controls(7,-0.7)and (10,-0.7)..(10.,-1.);
				\draw [line width=1.2pt,color=black] (7.,-2.)..controls(7,-1.7)and (9,-1.7)..(9.,-2.);
				\draw [line width=1.2pt,color=black] (7.,-2.)..controls(7,-1.7)and (10,-1.7)..(10.,-2.);
				\draw [line width=1.2pt,color=black] (7.,-3.)..controls(7,-2.7)and (9,-2.7)..(9.,-3.);
				\draw [line width=1.2pt,color=black] (7.,-3.)..controls(7.,-2.7)and (10.,-2.7)..(10.,-3.);

				\draw [line width=1.2pt,color=black] (-3.,2.)..controls(-3.3,2.)and (-3.3,4.)..(-3.,4.);
				\draw [line width=1.2pt,color=black] (-3.,2.)..controls(-3.3,2.)and (-3.3,5.)..(-3.,5.);
				\draw [line width=1.2pt,color=black] (-2.,2.)..controls(-2.3,2.)and (-2.3,4.)..(-2.,4.);
				\draw [line width=1.2pt,color=black] (-2.,2.)..controls(-2.3,2.)and (-2.3,5.)..(-2.,5.);
				\draw [line width=1.2pt,color=black] (-1.,2.)..controls(-1.3,2.)and (-1.3,4.)..(-1.,4.);
				\draw [line width=1.2pt,color=black] (-1.,2.)..controls(-1.3,2.)and (-1.3,5.)..(-1.,5.);
				\draw [line width=1.2pt,color=black] (0.,2.)..controls(-0.3,2.)and (-0.3,4.)..(0.,4.);
				\draw [line width=1.2pt,color=black] (0.,2.)..controls(-0.3,2.)and (-0.3,5.)..(0.,5.);
				\draw [line width=1.2pt,color=black] (2.,2.)..controls(1.7,2.)and (1.7,4.)..(2.,4.);
				\draw [line width=1.2pt,color=black] (2.,2.)..controls(1.7,2.)and (1.7,5.)..(2.,5.);
				\draw [line width=1.2pt,color=black] (3.,2.)..controls(2.7,2.)and (2.7,4.)..(3.,4.);
				\draw [line width=1.2pt,color=black] (3.,2.)..controls(2.7,2.)and (2.7,5.)..(3.,5.);
				\draw [line width=1.2pt,color=black] (4.,2.)..controls(3.7,2.)and (3.7,4.)..(4.,4.);
				\draw [line width=1.2pt,color=black] (4.,2.)..controls(3.7,2.)and (3.7,5.)..(4.,5.);
				\draw [line width=1.2pt,color=black] (5.,2.)..controls(4.7,2.)and (4.7,4.)..(5.,4.);
				\draw [line width=1.2pt,color=black] (5.,2.)..controls(4.7,2.)and (4.7,5.)..(5.,5.);
				\draw [line width=1.2pt,color=black] (7.,2.)..controls(6.7,2.)and (6.7,4.)..(7.,4.);
				\draw [line width=1.2pt,color=black] (7.,2.)..controls(6.7,2.)and (6.7,5.)..(7.,5.);
				\draw [line width=1.2pt,color=black] (8.,2.)..controls(7.7,2.)and (7.7,4.)..(8.,4.);
				\draw [line width=1.2pt,color=black] (8.,2.)..controls(7.7,2.)and (7.7,5.)..(8.,5.);
				\draw [line width=1.2pt,color=black] (9.,2.)..controls(8.7,2.)and (8.7,5.)..(9.,5.);
				\draw [line width=1.2pt,color=black] (9.,2.)..controls(8.7,2.)and (8.7,4.)..(9.,4.);
				\draw [line width=1.2pt,color=black] (10.,2.)..controls(9.7,2.)and (9.7,4.)..(10.,4.);
				\draw [line width=1.2pt,color=black] (10.,2.)..controls(9.7,2.)and (9.7,5.)..(10.,5.);
				\draw [line width=1.2pt,color=black] (-3.,-3.)..controls(-3.3,-3.)and (-3.3,-1.)..(-3.,-1.);
				\draw [line width=1.2pt,color=black] (-3.,-3.)..controls(-3.3,-3)and (-3.3,0.)..(-3.,0.);
				\draw [line width=1.2pt,color=black] (-2.,-3.)..controls(-2.3,-3.)and (-2.3,-1.)..(-2.,-1.);
				\draw [line width=1.2pt,color=black] (-2.,-3.)..controls(-2.3,-3.)and (-2.3,0.)..(-2.,0.);
				\draw [line width=1.2pt,color=black] (-1.,-3.)..controls(-1.3,-3.)and (-1.3,-1.)..(-1.,-1.);
				\draw [line width=1.2pt,color=black] (-1.,-3.)..controls(-1.3,-3.)and (-1.3,0.)..(-1.,0.);
				\draw [line width=1.2pt,color=black] (0.,-3.)..controls(-0.3,-3.)and (-0.3,-1.)..(0.,-1.);
				\draw [line width=1.2pt,color=black] (0.,-3.)..controls(-0.3,-3.)and (-0.3,0.)..(0.,0.);
				\draw [line width=1.2pt,color=black] (2.,-3.)..controls(1.7,-3.)and (1.7,-1.)..(2.,-1.);
				\draw [line width=1.2pt,color=black] (2.,-3.)..controls(1.7,-3.)and (1.7,0.)..(2,0.);
				\draw [line width=1.2pt,color=black] (3.,-3.)..controls(2.7,-3.)and (2.7,-1.)..(3.,-1.);
				\draw [line width=1.2pt,color=black] (3.,-3.)..controls(2.7,-3.)and (2.7,0.)..(3.,0.);
				\draw [line width=1.2pt,color=black] (4.,-3.)..controls(3.7,-3.)and (3.7,-1.)..(4.,-1.);
				\draw [line width=1.2pt,color=black] (4.,-3.)..controls(3.7,-3.)and (3.7,0.)..(4.,0.);
				\draw [line width=1.2pt,color=black] (5.,-3.)..controls(4.7,-3.)and (4.7,-1.)..(5.,-1.);
				\draw [line width=1.2pt,color=black] (5.,-3.)..controls(4.7,-3.)and (4.7,0.)..(5.,0.);
				\draw [line width=1.2pt,color=black] (7.,-3.)..controls(6.7,-3.)and (6.7,-1.)..(7.,-1.);
				\draw [line width=1.2pt,color=black] (7.,-3.)..controls(6.7,-3.)and (6.7,0.)..(7.,0.);
				\draw [line width=1.2pt,color=black] (8.,-3.)..controls(7.7,-3.)and (7.7,-1.)..(8.,-1.);
				\draw [line width=1.2pt,color=black] (8.,-3.)..controls(7.7,-3.)and (7.7,0.)..(8.,0.);
				\draw [line width=1.2pt,color=black] (9.,-3.)..controls(8.7,-3.)and (8.7,-1.)..(9.,-1.);
				\draw [line width=1.2pt,color=black] (9.,-3.)..controls(8.7,-3.)and (8.7,0.)..(9.,0.);
				\draw [line width=1.2pt,color=black] (10.,-3.)..controls(9.7,-3.)and (9.7,-1.)..(10.,-1.);
				\draw [line width=1.2pt,color=black] (10.,-3.)..controls(9.7,-3.)and (9.7,0.)..(10.,0.);
				\begin{scriptsize}
					\draw [fill=white] (-3.,5.) circle (2.5pt);
					\draw [fill=white] (-3.,4.) circle (2.5pt);
					\draw [fill=white] (-3.,3.) circle (2.5pt);
					\draw [fill=white] (-3.,2.) circle (2.5pt);
					\draw [fill=red] (-2.,2.) circle (2.5pt);
					\draw (-2.25,1.95) node[anchor=north west] {\footnotesize$d_1$};
					\draw [fill=white] (-2.,3.) circle (2.5pt);
					\draw [fill=white] (-2.,4.) circle (2.5pt);
					\draw [fill=white] (-2.,5.) circle (2.5pt);
					\draw [fill=white] (-1.,5.) circle (2.5pt);
					\draw [fill=white] (-1.,4.) circle (2.5pt);
					\draw [fill=white] (-1.,3.) circle (2.5pt);
					\draw [fill=blue] (-1.,2.) circle (2.5pt);
					\draw (-1.25,1.85) node[anchor=north west] {\footnotesize$s_1$};
					\draw [fill=blue] (0.,2.) circle (2.5pt);
					\draw (-0.25,1.85) node[anchor=north west] {\footnotesize$s_2$};
					\draw [fill=white] (0.,3.) circle (2.5pt);
					\draw [fill=white] (0.,4.) circle (2.5pt);
					\draw [fill=white] (0.,5.) circle (2.5pt);
					\draw [fill=white] (2.,2.) circle (2.5pt);
					\draw [fill=white] (2.,3.) circle (2.5pt);
					\draw [fill=white] (2.,4.) circle (2.5pt);
					\draw [fill=white] (2.,5.) circle (2.5pt);
					\draw [fill=red] (3.,2.) circle (2.5pt);
					\draw (2.75,1.95) node[anchor=north west] {\footnotesize$d_1$};
					\draw [fill=white] (3.,3.) circle (2.5pt);
					\draw [fill=white] (3.,4.) circle (2.5pt);
					\draw [fill=white] (3.,5.) circle (2.5pt);
					\draw [fill=blue] (4.,2.) circle (2.5pt);
					\draw (4.1,2.15) node[anchor=north west] {\footnotesize$s_1$};
					\draw [fill=blue] (4.,3.) circle (2.5pt);
					\draw (4.1,3.15) node[anchor=north west] {\footnotesize$s_2$};
					\draw [fill=white] (4.,4.) circle (2.5pt);
					\draw [fill=white] (4.,5.) circle (2.5pt);
					\draw [fill=red] (5.,2.) circle (2.5pt);
					\draw (4.75,1.95) node[anchor=north west] {\footnotesize$d_2$};
					\draw [fill=white] (5.,3.) circle (2.5pt);
					\draw [fill=white] (5.,4.) circle (2.5pt);
					\draw [fill=white] (5.,5.) circle (2.5pt);
					\draw [fill=white] (7.,2.) circle (2.5pt);
					\draw [fill=red] (7.,3.) circle (2.5pt);
					\draw (6.25,3.25) node[anchor=north west] {\footnotesize$d_1$};
					\draw [fill=blue] (7.,4.) circle (2.5pt);
					\draw (6.25,4.25) node[anchor=north west] {\footnotesize$s_1$};
					\draw [fill=blue] (7.,5.) circle (2.5pt);
					\draw (6.25,5.25) node[anchor=north west] {\footnotesize$s_2$};
					\draw [fill=white] (8.,2.) circle (2.5pt);
					\draw [fill=white] (8.,3.) circle (2.5pt);
					\draw [fill=white] (8.,4.) circle (2.5pt);
					\draw [fill=white] (8.,5.) circle (2.5pt);
					\draw [fill=white] (9.,2.) circle (2.5pt);
					\draw [fill=white] (9.,3.) circle (2.5pt);
					\draw [fill=white] (9.,4.) circle (2.5pt);
					\draw [fill=white] (9.,5.) circle (2.5pt);
					\draw [fill=white] (10.,2.) circle (2.5pt);
					\draw [fill=white] (10.,3.) circle (2.5pt);
					\draw [fill=white] (10.,4.) circle (2.5pt);
					\draw [fill=white] (10.,5.) circle (2.5pt);
					\draw [fill=red] (-3.,0.) circle (2.5pt);
					\draw (-3.75,0.3) node[anchor=north west] {\footnotesize$d_2$};
					\draw [fill=blue] (-3.,-1.) circle (2.5pt);
					\draw (-3.75,-0.7) node[anchor=north west] {\footnotesize$s_1$};
					\draw [fill=red] (-3.,-2.) circle (2.5pt);
					\draw (-3.75,-1.7) node[anchor=north west] {\footnotesize$d_1$};
					\draw [fill=white] (-3.,-3.) circle (2.5pt);
					\draw [fill=white] (-2.,0.) circle (2.5pt);
					\draw [fill=blue] (-2.,-1.) circle (2.5pt);
					\draw (-1.85,-0.85) node[anchor=north west] {\footnotesize$s_2$};
					\draw [fill=white] (-2.,-2.) circle (2.5pt);
					\draw [fill=white] (-2.,-3.) circle (2.5pt);
					\draw [fill=white] (-1.,0.) circle (2.5pt);
					\draw [fill=white] (-1.,-1.) circle (2.5pt);
					\draw [fill=white] (-1.,-2.) circle (2.5pt);
					\draw [fill=white] (-1.,-3.) circle (2.5pt);
					\draw [fill=white] (0.,0.) circle (2.5pt);
					\draw [fill=white] (0.,-1.) circle (2.5pt);
					\draw [fill=white] (0.,-2.) circle (2.5pt);
					\draw [fill=white] (0.,-3.) circle (2.5pt);
					\draw [fill=white] (2.,0.) circle (2.5pt);
					\draw [fill=blue] (2.,-1.) circle (2.5pt);
					\draw (1.25,-0.75) node[anchor=north west] {\footnotesize$s_1$};
					\draw [fill=white] (2.,-2.) circle (2.5pt);
					\draw [fill=white] (2.,-3.) circle (2.5pt);
					\draw [fill=white] (3.,0.) circle (2.5pt);
					\draw [fill=white] (3.,-1.) circle (2.5pt);
					\draw [fill=red] (3.,-2.) circle (2.5pt);
					\draw (3.15,-1.8) node[anchor=north west] {\footnotesize$d_1$};
					\draw [fill=white] (3.,-3.) circle (2.5pt);
					\draw [fill=white] (4.,0.) circle (2.5pt);
					\draw [fill=white] (4.,-1.) circle (2.5pt);
					\draw [fill=white] (4.,-2.) circle (2.5pt);
					\draw [fill=white] (4.,-3.) circle (2.5pt);
					\draw [fill=white] (5.,0.) circle (2.5pt);
					\draw [fill=white] (5.,-1.) circle (2.5pt);
					\draw [fill=white] (5.,-2.) circle (2.5pt);
					\draw [fill=white] (5.,-3.) circle (2.5pt);
					\draw [fill=white] (7.,0.) circle (2.5pt);
					\draw [fill=white] (7.,-1.) circle (2.5pt);
					\draw [fill=white] (7.,-2.) circle (2.5pt);
					\draw [fill=red] (7.,-3.) circle (2.5pt);
					\draw (6.75,-3.1) node[anchor=north west] {\footnotesize$d_1$};
					\draw [fill=white] (8.,0.) circle (2.5pt);
					\draw [fill=white] (8.,-1.) circle (2.5pt);
					\draw [fill=white] (8.,-2.) circle (2.5pt);
					\draw [fill=blue] (8.,-3.) circle (2.5pt);
					\draw (7.75,-3.2) node[anchor=north west] {\footnotesize$s_1$};
					\draw [fill=white] (9.,0.) circle (2.5pt);
					\draw [fill=white] (9.,-1.) circle (2.5pt);
					\draw [fill=white] (9.,-2.) circle (2.5pt);
					\draw [fill=white] (9.,-3.) circle (2.5pt);
					\draw [fill=white] (10.,0.) circle (2.5pt);
					\draw [fill=white] (10.,-1.) circle (2.5pt);
					\draw [fill=white] (10.,-2.) circle (2.5pt);
					\draw [fill=white] (10.,-3.) circle (2.5pt);
					\draw (-2,1.5) node[anchor=north west] {\normalsize$(a)$};
					\draw (3,1.5) node[anchor=north west] {\normalsize$(b)$};
					\draw (8,1.5) node[anchor=north west] {\normalsize$(c)$};
					\draw (-2,-3.5) node[anchor=north west] {\normalsize$(d)$};
					\draw (3,-3.5) node[anchor=north west] {\normalsize$(e)$};
					\draw (8,-3.5) node[anchor=north west] {\normalsize$(f)$};
				\end{scriptsize}
			\end{tikzpicture}
			\caption{An illustration for Staller's winning strategy in the proof proof of Theorem \ref{prop:out_star} Case 3. Blue and red vertices denote the moves of Staller and Dominator, respectively.}
			\label{fig:StallStar Dgame}
		\end{center}
	\end{figure}
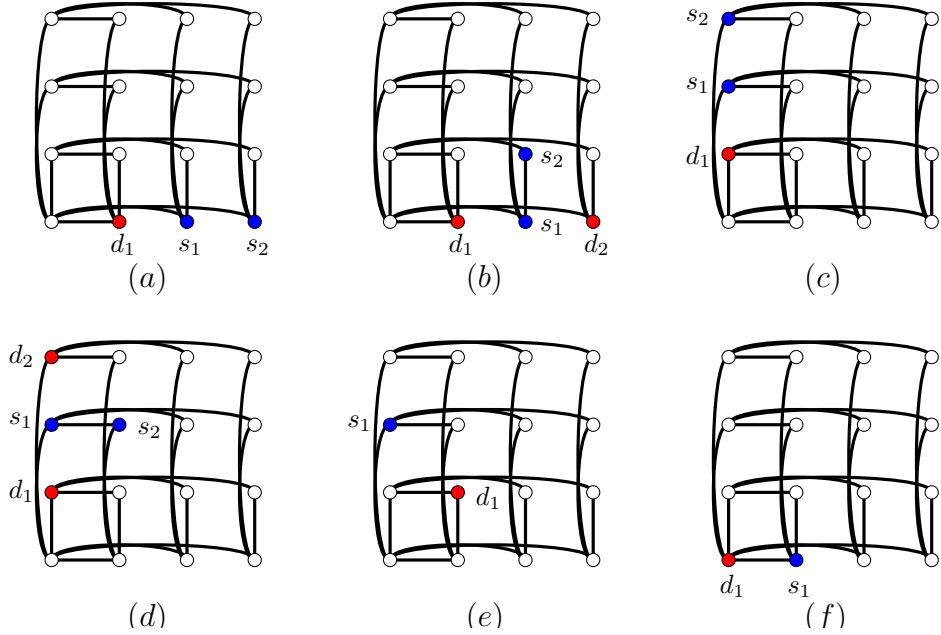

	\textbf{Case 3.1.} $d_1=(a,1)$. Then there are two unplayed layers $(K_{1,3})^2$ and $(K_{1,3})^3$. Then Staller plays $s_1=(a,2)$. 	
	See Figure \ref{fig:StallStar Dgame}\,(a)-(b).
	\begin{itemize}
		\item If $d_2 \ne (a,3)$, then Staller replies $s_2=(a,3)$. After Dominator plays $d_3$, there exists $i\in[3]$ such that $(i,b), (i,2), (i,3)$ are unplayed. Then Staller plays $(i,b)$, and she will win in her next move by claiming the closed neighborhood of $(i,2)$ or $(i,3)$.
		
		\item If $d_2=(a,3)$, then Staller replies $s_2=(1,2)$. It forces Dominator to play $d_3=(1,b)$. Then Staller can play all neighbors of $(a,2)$, and she will win the game.
		
	\end{itemize}
	
	\textbf{Case 3.2.} $d_1 = (1,b)$. Then Staller replies $s_1=(2,b)$. By symmetry and commutativity, we can see $(1,b)$ as $(a,1)$ and $(2,b)$ as $(2,a)$. Then we can apply the above strategy which implies the same outcome. See Figure \ref{fig:StallStar Dgame}\,(c)-(d).
	
	\textbf{Case 3.3.} $d_1=(1,1)$. Then Staller plays $s_1=(2,b)$. The continuation of the game Staller applies the same strategy as Case 3.1.
	See Figure \ref{fig:StallStar Dgame}\,(e).
	
	\textbf{Case 3.4.} $d_1=(a,b)$.  There are three unplayed layers $(K_{1,3})^1$, $(K_{1,3})^2$, and $(K_{1,3})^3$. Then Staller plays $s_1=(a,1)$ and plays $s_2=(a,j)$, where $j \in \{2,3\}$. After the third move of Dominator, there exists $i\in[3]$ such that $(i,b), (i,1), (i,j)$ are unplayed. Thus Staller plays $(i,b)$ and she will win in her next move by claiming the closed neighborhood of $(i,1)$ or $(i,j)$.
	See Figure \ref{fig:StallStar Dgame}\,(f).
	
	Thus Staller has a strategy to win the game on $Z$ which implies that $o(Z) = \cal S$.
	\paragraph{Case 4.} $m \ge 3,\ n \ge 4$.  By Corollary \ref{cor:skip}\,(ii), to show that Staller has a winning strategy in the D-game on $Z$ that implies $o(Z) = \cal S$. Since $m \ge 3$, after the first move of Dominator, there are two unplayed layers $^i(K_{1,n})$ and $^{i'}(K_{1,n})$, where $i,i' \in [m]$. Then Staller plays $(i,b)$ in her first move.
	Since $n \ge 4$, after the second move of Dominator, there are two unplayed layers $(K_{1,m})^j$ and $(K_{1,m})^{j'}$, where $j,j' \in [n]$.
	Then Staller replies $(a,j)$ in her second move and it makes Dominator play $(i,j)$, otherwise, Staller will win in her next move. By this strategy, $(K_{1,3})^{j'}$ is still unplayed at this moment, so, Staller plays $(a,j')$. If Dominator does not reply by playing $(i,j')$, then Staller will win the game by playing $(i,j')$. If Dominator replies $(i,j')$, then Staller responds at $(i',b)$ and she will win the game in her next turn. 
	
	\medskip
	We conclude that $o(K_{1,m} \cp K_{1,n}) = \cS$, where $m \ge 3, \, n\ge 4.$
	\qed
	\medskip
	
	\begin{theorem}
		If  $n \ge 3$, then
		\begin{align*}
			\gmb (P_3 \cp K_{1,n}) = \gsmb' (P_3 \cp K_{1,n}) = n+2.
		\end{align*}
	\end{theorem}
	\proof  Assume that $n \ge 3$ and $Z = P_3 \cp K_{1,n}$. Recall that $a, b$ are the central vertices of $P_3$ and $K_{1,n}$, respectively.
	
	First, we consider the D-game. We will show that  $\gmb (Z) \le n+2$.
	Dominator starts the game by playing $(1,b)$. Then there is a matching $M$ in $Z-\{(1,b)\}$ such that $V(Z)\setminus V(M) \subseteq N_Z[(1,b)]$ and $|M| = n+1$. By Lemma~\ref{lem:pairing} and Remark~\ref{re:pairing}\,, Dominator has a strategy to win the game within $n+2$ moves. Thus $\gmb (Z) \le n+2$.
	
	Now we will show that Staller has a strategy ensure that Dominator needs to play at least $n+2$ moves.
	By symmetry, we consider the following cases.
	\begin{itemize}
		\item If $d_1=(1,b)$, then Staller plays $s_1=(2,1)$.
		\item If $d_1= (a,b)$, then Staller plays $s_1=(a,1)$.
		\item If $d_1=(1,1)$, then Staller plays  $s_1=(2,b)$.
		\item If $d_1=(a,1)$, then Staller plays $s_1=(a,b)$.
	\end{itemize}
	In the second turn, Staller selects an unplayed vertex $s_2$ in the layer $(P_3)^b$. For the continuation of the game, she will try to play a vertex $(a,j)$ for $j \in [n]$ if it is possible, otherwise she plays an arbitrary unplayed vertex.
	One can see that Dominator needs at least two moves to dominate the layers $(P_3)^b$ and $(P_3)^1$, and $n-1$ moves to dominate the rest of the graph. Thus $\gmb (Z) \ge n+2$ and hence $\gmb(Z) = n+2$.
	
	Next, we consider the S-game. By the proof of Theorem~\ref{prop:out_star} Case 2,  Staller has a strategy to win the game within $n+2$ moves.
	
	It remains to show that $\gsmb' (Z) \ge n+2$. Notice that closed neighborhoods in $Z$ are of size $3$, $4$, $n+2$, and $n+3$. 
	By symmetry, it is enough  to consider the following cases.
	If $s'_1 \in \{(a,b), (a,1), (1,1), (2,b)\}$, then Dominator plays $d'_1=(1,b)$. If $s'_1 \in \{(1,b), (2,1)\}$, then Dominator plays $d'_1=(2,b)$. 
	
	\begin{itemize}

		\item If $(2,b)$ is unplayed after the second move of Staller, then Dominator plays $d'_2=(2,b)$. For the continuation of the game, when Staller plays $(i,j)$, where $i \in [2], j\in [n]$ he replies by playing $(a,j)$ if it is possible. Otherwise, he plays an arbitrary unplayed vertex. If Staller plays $(a, j)$, where $j \in [n]$, he replies by playing $(i, j)$, where $i \in [2]$ if it is possible, otherwise he plays an arbitrary unplayed vertex.
		
		\item Assume that Staller occupies $(2,b)$. If $(a,j)$, where $j \in [n]$ was played by Staller, then Dominator plays $d'_2=(2,j)$. Otherwise, he plays $(a,j)$ such that $(i,j)$, where $i\in [2]$ was played by Staller. For the continuation of the game, when Staller plays $(i',j')$, where $i'\in [2], j'\in [n]$ he replies by playing $(a,j')$ if it is possible. Otherwise, he plays an arbitrary unplayed vertex. When Staller plays $(a, j')$, where $j' \in [n]$ and $(a,b)$ is unplayed, then he replies $(a,b)$. If $(a,b)$ is played by Staller, then he replies by playing $(1, j')$. Otherwise he plays an arbitrary unplayed vertex.
	\end{itemize}
	
	By above strategy, Staller cannot claim neither a closed neighborhood of size $3$ nor $4$. Thus Staller needs to play at least $n+2$ moves. 
	
	\qed
	\medskip
	
	\begin{theorem}
		If $m \ge n \ge3$, then
		\begin{align*}
			\gsmb (K_{1,m} \cp K_{1,n}) &= \begin{cases}
				5; & n=3,\\
				4; & n \ge 4,
			\end{cases}
		\end{align*}
		
		\noindent and  
		\begin{align*}
			\gsmb' (K_{1,m} \cp K_{1,n}) = 4.
		\end{align*}
	\end{theorem}
	\proof Assume that $ m \ge n \ge3$ and $Z=K_{1,m} \cp K_{1,n}$.  Recall that $a, b$ are the central vertices of $K_{1,m}$ and $K_{1,n}$, respectively. Observe that every vertex is in the closed neighborhood of size $3, m+2, n+2$, and $m+n+1$, respectively. We first consider the D-game.
	\paragraph{Case 1.} $n=3$ and now $Z = K_{1,m} \cp K_{1,3}$.
	By the proof of Theorem~\ref{prop:out_star} Case 3 and Case 4,  Staller has a strategy to win the D-game within five moves. It implies that $\gsmb(Z) \le 5$. It remains to show that Staller needs to play at least five moves to win the game. Assume Dominator plays $d_1=(a,1)$ in his first move. Assume that Staller plays $s_1$ in her first move.
	\begin{itemize}
		\item $s_1= (a,2)$. Then Dominator replies by playing $d_2=(a,3)$. For each turn, if Staller plays $(i,j)$, then Dominator replies by selecting $(i, b)$ if it is possible, otherwise he plays $(i, 2)$ if it is available or any arbitrary vertex if $(i,2)$ was played earlier. Thus Staller cannot claim a closed neighborhood of size $3$. Hence, Staller needs to play at least five moves.
		
		\item $s_1\ne (a,2)$. Then Dominator plays $d_2=(a,2)$. After the second move of Staller, if $(a,3)$ is unplayed, then Dominator will play $d_3=(a,3)$. 
		Assume that Staller already played $(a,3)$. 
		If $(i,b)$, where $i \in [3]$ was played by Staller, then he will play $d_3=(i,3)$. If  $(i,3)$, where $i \in [3]$ was played by Staller, then he will play $d_3=(i,b)$.
		Otherwise, Dominator will play an arbitrary unplayed vertex. 
		Then he applies this strategy to the continuation of the game. 
		Thus Staller cannot win the game by claiming a closed neighborhood of size $3$. So, she needs to play at least five moves.
		
	\end{itemize}
	Therefore, $\gsmb (Z) \ge 5$. We conclude that $\gsmb (Z) = 5.$
	
	\paragraph{Case 2.} $n \ge 4$ and now $Z =K_{1,m} \cp K_{1,n}$. By Lemma~\ref{cor:skip}, $\gsmb(Z) \ge \gsmb' (Z)$.
	It remains to show that $\gsmb (Z)  \le 4$ and $ \gsmb' (Z) \ge 4$. We first provide a strategy for Staller to win the D-game in four moves. After the first move of Dominator, there exists $j \in [3]$ such that $(a, j)$ is unplayed. Then Staller plays $(a,j)$. Similarly, after the second move of Dominator, there exists $j' \in [3]$ such that $j' \ne j$ and $(a,j')$ is unplayed. Then Staller plays  $(a,j')$. Since $m \ge n \ge 4$, after the third move of Dominator, there is $(i,b)$ such that $(i,b), (i,j), (i,j')$ are unplayed. So, Staller plays $(i,b)$ and she can win the game in her next move. Thus $\gsmb (Z)  \le 4$.
	
	For the S-game, we will show that Staller cannot win the game in three moves. Observe that every vertex is in a closed neighborhood of size $3$ except $(a,b)$. By symmetry, we consider the following cases.
	\begin{itemize}
		\item If $s_1'=(a,b)$, then she needs to play at least four moves to claim a closed neighborhood.
		\item If $s_1'\in \{(a,1),(1,b)\}$, then Dominator replies by playing $(1,1)$. Thus Staller needs to play at least three more moves to claim a closed neighborhood.
		\item If $s_1'=(1,1)$, then Dominator replies by playing $(1,a)$. Thus Staller needs to play at least three more moves to claim a closed neighborhood.
	\end{itemize}
	It implies that $ \gsmb'(Z) \ge 4$.
	Therefore, $ \gsmb'(Z) = \gsmb(Z) = 4$.
	\qed
	
	\section*{Acknowledgements}
	The author would like to express our sincere gratitude to Prof.\ Sandi Klav\v{z}ar and Assoc.\ Prof.\ Csilla Bujt\'{a}s, for their invaluable guidance and support.

\end{document}